\documentclass[11pt,a4paper]{article}

\usepackage{amsthm,amsmath,amsfonts,amssymb}
\usepackage{mathtools,mathrsfs,enumitem}
\usepackage[numbers]{natbib}
\usepackage[hidelinks]{hyperref}

\usepackage[margin=1in]{geometry}

\theoremstyle{definition}
\newtheorem{definition}{Definition}[section]
\newtheorem{assumption}{Assumption}[section]
\newtheorem{example}{Example}[section]

\theoremstyle{plain}
\newtheorem{theorem}{Theorem}[section]
\newtheorem{lemma}{Lemma}[section]
\newtheorem{proposition}{Proposition}[section]
\newtheorem{corollary}{Corollary}[section]

\theoremstyle{definition}
\newtheorem{remark}{Remark}[section]

\newcommand{\R}{\mathbb{R}}
\newcommand{\Z}{\mathbb{Z}}
\newcommand{\N}{\mathbb{N}}

\newcommand{\T}{\mathbb{T}}

\newcommand{\E}{\mathbb{E}}
\newcommand{\Var}{\operatorname{Var}}
\newcommand{\Cov}{\operatorname{Cov}}

\newcommand{\1}{\mathbf{1}}
\newcommand{\dd}{\,\mathrm{d}}
\DeclareMathOperator{\supp}{supp}

\title{Central Limit Theorems for Iterated Convolution Fields Built on a Poisson Point Process}
\author{%
  K. V. Hlyniana\textsuperscript{1,2} \and
  A. A. Dorogovtsev\textsuperscript{2}}
\date{%
  \small\textsuperscript{1}School of Mathematics, Jilin University\\
  \small\textsuperscript{2}Institute of Mathematics, National Academy of Sciences of Ukraine\\[4pt]
  \small\href{mailto:glinkate@gmail.com}{glinkate@gmail.com} (K. V. Hlyniana)\\
  \small\href{mailto:andrey.dorogtsev@gmail.com}{andrey.dorogtsev@gmail.com} (A. A. Dorogovtsev)}

\begin{document}

\maketitle

\begin{abstract}
We study recursive convolution fields built over a homogeneous Poisson point process on the real line. The construction includes both linear and nonlinear recursive weights generated by a compactly supported interaction kernel. For bounded periodic test functions, we investigate additive statistics over growing observation windows. We prove central limit theorems at fixed recursive depth and when the depth grows with the window size, identify the covariance of the fixed-depth Gaussian limit, and describe its Fourier representation and growth with the recursive depth.
\end{abstract}



\section{Introduction}
\label{sec:introduction}
%===========================================================

In this paper we study random fields generated by repeatedly applying a
local convolution mechanism to one fixed Poisson point configuration. Our
main objective is to describe the macroscopic fluctuations of the resulting
weighted fields when the observation window becomes large. We consider both
a linear propagation rule and a nonlinear reinforcement rule, and we study
two asymptotic regimes: fixed recursive depth and recursive depth increasing
with the observation scale.

The objects of interest are linear statistics obtained by integrating a
deterministic test function against the recursively weighted point
configuration. At fixed depth, we prove joint central limit theorems for
arbitrary finite collections of such statistics and identify the covariance
of the limiting Gaussian linear family. This covariance has a direct
statistical interpretation. Once the integrals against finitely many basis
functions are known, the Gaussian limit permits the integral against another
test function to be approximated by its projection onto the observed
coordinates. In particular, Fourier diagonalization yields an explicit
decomposition into an observable component and an independent residual
component, while the omitted covariance eigenvalues quantify the
mean-square approximation error. Thus the limiting covariance describes not
only the fluctuations of the recursive field, but also how accurately one
linear functional can be recovered from a finite collection of other
linear functionals.

We now introduce the two recursive models. Let \(\mu\) be a homogeneous
Poisson point process on \(\R\) with intensity \(\lambda>0\). Let
\(\varphi:\R\to[0,\infty)\) be measurable and bounded, and suppose that
\begin{equation}
\label{eq:intro-kernel-support}
\supp\varphi\subset[-R,R]
\end{equation}
for some \(R<\infty\). We also assume that
\begin{equation}
\label{eq:intro-alpha-zero}
\alpha_0:=\varphi(0)>0.
\end{equation}
For a nonnegative measurable field \(H\), define the random convolution
operator
\begin{equation}
\label{eq:intro-random-convolution-operator}
(\mathcal T_\mu H)(v):=\int_{\R}\varphi(u-v)H(u)\,\mu(\dd u).
\end{equation}

The linear recursion is
\begin{equation}
\label{eq:intro-linear-recursion}
F_0\equiv1,\qquad F_{n+1}=\mathcal T_\mu F_n,\qquad n\ge0,
\end{equation}
or equivalently,
\[
F_{n+1}(v)=\int_{\R}\varphi(u-v)F_n(u)\,\mu(\dd u).
\]
Conditionally on \(\mu\), regard the Poisson points as sites of a weighted
spatial network and assign the directed interaction weight
\[
K(v,u):=\varphi(u-v)
\]
to the ordered pair \((v,u)\). Iterating
\eqref{eq:intro-linear-recursion} gives
\begin{equation}
\label{eq:intro-linear-walk-expansion}
F_n(v)=\int_{\R^n}\prod_{j=1}^{n}\varphi(u_j-u_{j-1})\,
\mu(\dd u_1)\cdots\mu(\dd u_n),\qquad u_0:=v.
\end{equation}
Thus, when \(v\in\mu\), the field \(F_n(v)\) is the total weight of all
ordered \(n\)-step walks rooted at \(v\). Repeated sites and self-loops are
allowed, with the latter contributing through
\(\alpha_0=\varphi(0)\). The linear recursion therefore represents repeated
propagation and aggregation through a single random interaction network.

The nonlinear recursion is
\begin{equation}
\label{eq:intro-nonlinear-recursion}
A_0\equiv1,\qquad A_{n+1}=A_n\,\mathcal T_\mu A_n,\qquad n\ge0,
\end{equation}
that is,
\[
A_{n+1}(v)=A_n(v)\int_{\R}\varphi(u-v)A_n(u)\,\mu(\dd u).
\]
Define the weighted random measure
\begin{equation}
\label{eq:intro-eta-n}
\eta_n(\dd v):=A_n(v)\mu(\dd v).
\end{equation}
Its support remains equal to the original Poisson configuration, while the
mass assigned to each point is updated recursively. Since
\(\eta_0=\mu\), equation \eqref{eq:intro-nonlinear-recursion} is equivalent
to
\begin{equation}
\label{eq:intro-eta-recursion}
\eta_{n+1}(\dd v)=\left(\int_{\R}\varphi(u-v)\eta_n(\dd u)\right)\eta_n(\dd v).
\end{equation}
The convolution in \eqref{eq:intro-eta-recursion} represents the weighted
influence of the current measure near \(v\). The new mass at \(v\) is the
product of the mass already present there and this local influence. The
nonlinear dynamics therefore produces multiplicative local reinforcement
without creating new sites or moving the original points.

The two models have the same first iteration but generate substantially
different combinatorial complexity at larger depths. Repeated substitution
represents both \(F_n(v)\) and \(A_n(v)\) as multiple integrals with respect
to \(\mu\). Let \(d_n\) denote the number of integration variables in the
corresponding representation. Then
\begin{equation}
\label{eq:intro-degree}
d_n=
\begin{cases}
n, & \text{for the linear recursion},\\[2mm]
2^n-1, & \text{for the nonlinear recursion}.
\end{cases}
\end{equation}
Each linear iteration introduces one new integration variable. In the
nonlinear recursion, the two factors of multiplicity \(d_n\) are combined
with one additional integration variable, so that
\(d_{n+1}=2d_n+1\).

We use the unified notation
\begin{equation}
\label{eq:intro-Wn}
W_n:=
\begin{cases}
F_n, & \text{for the linear recursion},\\
A_n, & \text{for the nonlinear recursion}.
\end{cases}
\end{equation}
For a bounded measurable \(1\)-periodic function \(f:\R\to\R\), let
\[
\Lambda_N:=[-N,N],\qquad |\Lambda_N|=2N,
\]
and define
\begin{equation}
\label{eq:intro-QnN}
Q_{n,N}(f):=\int_{\Lambda_N}f(v)W_n(v)\,\mu(\dd v).
\end{equation}
Thus \(Q_{n,N}(f)\) tests the weighted point measure
\(W_n(v)\mu(\dd v)\) against the observation profile \(f\). In the
nonlinear model,
\[
Q_{n,N}(f)=\int_{\Lambda_N}f(v)\,\eta_n(\dd v).
\]

The finite support of the kernel implies that \(W_n(v)\) depends only on
the restriction of \(\mu\) to
\begin{equation}
\label{eq:intro-locality-region}
[v-nR,v+nR].
\end{equation}
Hence the localization radius grows linearly with the recursive depth. For
fixed \(n\), the centered statistic \(Q_{n,N}(f)-\E Q_{n,N}(f)\) can be
decomposed into a stationary finitely dependent sequence of block
contributions. If \(n=n_N\) increases with \(N\), the same decomposition
produces a triangular array whose dependence range also increases with
\(N\).

The present problem is related to several established directions in the
limit theory of Poisson functionals. Convolutions of deterministic kernels
with Poisson processes produce classical shot-noise fields, while local
scores and multiple interactions give rise to random measures, geometric
statistics, and Poisson \(U\)-statistics. Palm--Mecke calculus, Poisson
integration, and chaos representations are developed systematically by Last
and Penrose \cite{LastPenrose2011,LastPenrose2017}.

Gaussian approximation on Poisson space has been studied through
Malliavin--Stein methods, stabilization, and finite-chaos expansions.
Malliavin--Stein bounds for Poisson functionals were obtained by Peccati,
Sol\'e, Taqqu, and Utzet \cite{PeccatiSoleTaqquUtzet2010}, with
multidimensional extensions by Peccati and Zheng
\cite{PeccatiZheng2010}. Reitzner and Schulte
\cite{ReitznerSchulte2013} studied Poisson \(U\)-statistics and their finite
Wiener--It\^o chaos expansions, while Last, Penrose, Schulte, and Th\"ale
\cite{LastPenroseSchulteThaele2014} established multivariate central limit
theorems for finite-chaos Poisson functionals. Spatial central limit
theorems based on stabilization and local dependence were developed, among
others, by Penrose \cite{Penrose2005} and Last, Peccati, and Schulte
\cite{LastPeccatiSchulte2016}. Random geometric graphs and related spatial
networks are treated in Penrose \cite{Penrose2003}.

Related theories extend beyond strictly finite-range Poisson functionals.
Lachi\`eze-Rey \cite{LachiezeRey2019} studied nonlocal geometric
functionals and excursion statistics of Poisson shot-noise fields.
B{\l}aszczyszyn, Yogeshwaran, and Yukich
\cite{BlaszczyszynYogeshwaranYukich2019} obtained laws of large numbers,
variance asymptotics, and central limit theorems for stabilizing statistics
of stationary point processes with fast decay of correlations.

The main distinction of the present setting is that the local score is not
prescribed in advance. It is generated recursively from the same Poisson
environment, and the dependence created at one depth is propagated through
all subsequent iterations. Moreover, the recursive depth may increase
together with the observation window. Consequently, the algebraic order of
the functional, the localization radius, the relevant moments, and the
variance may all change simultaneously. The growing-depth problem therefore
requires estimates that remain effective while the statistic itself becomes
progressively more complex.

We now summarize the principal results. Introduce the depth-\(n\) weighted
point measure
\[
\nu_n(\dd v):=W_n(v)\mu(\dd v),
\]
so that
\[
Q_{n,N}(f)=\int_{\Lambda_N}f(v)\,\nu_n(\dd v).
\]
Its second-order structure separates into a diagonal contribution, arising
when the same Poisson point appears in both factors, and an off-diagonal
contribution, arising from two distinct points whose recursively generated
neighbourhoods overlap. Let \(\E^0\) and \(\E^{0,h}\) denote expectation
with respect to the one-point and two-point Palm distributions. For fixed
\(n\), define
\begin{equation}
\label{eq:intro-an}
a_n:=\lambda\E^0[W_n(0)^2]
\end{equation}
and, for \(h\neq0\),
\begin{equation}
\label{eq:intro-bn}
b_n(h):=\lambda^2\left(\E^{0,h}[W_n(0)W_n(h)]-\bigl(\E^0W_n(0)\bigr)^2\right),
\qquad b_n(0):=0.
\end{equation}
The coefficient \(a_n\) determines the diagonal covariance, whereas
\(b_n(h)\) describes the off-diagonal covariance density at separation
\(h\).

The first result establishes a joint Gaussian limit at every fixed depth and
identifies its covariance.

\begin{theorem}[Fixed-depth Gaussian limit]
\label{thm:intro-fixed-depth}
For every fixed \(n\ge0\), the function \(b_n\) is even, nonnegative, and
supported on \([-2nR,2nR]\). For bounded \(1\)-periodic functions \(f\) and
\(g\),
\begin{equation}
\label{eq:intro-Sigma-n}
\begin{aligned}
\Sigma_n(f,g)
&:=\lim_{N\to\infty}\frac{1}{|\Lambda_N|}
\Cov\bigl(Q_{n,N}(f),Q_{n,N}(g)\bigr)\\
&=a_n\int_0^1f(x)g(x)\,\dd x
+\int_0^1\int_{\R}f(x)g(x+h)b_n(h)\,\dd h\,\dd x.
\end{aligned}
\end{equation}
The form \(\Sigma_n\) extends to a bounded, symmetric, and nonnegative
bilinear form on \(L^2([0,1])\) and defines a centered Gaussian linear
family
\[
\{\xi_n(f):f\in L^2([0,1])\},
\qquad
\E[\xi_n(f)\xi_n(g)]=\Sigma_n(f,g).
\]
Moreover, for every \(m\ge1\) and every collection of bounded
\(1\)-periodic functions \(f_1,\ldots,f_m\),
\begin{equation}
\label{eq:intro-fixed-depth-fdd}
\left(
\frac{Q_{n,N}(f_j)-\E Q_{n,N}(f_j)}{\sqrt{|\Lambda_N|}}
\right)_{j=1}^{m}
\Longrightarrow
\bigl(\xi_n(f_j)\bigr)_{j=1}^{m},
\qquad N\to\infty.
\end{equation}
\end{theorem}

The recursive dependence therefore survives in the Gaussian limit through
\(a_n\) and \(b_n\). Formula \eqref{eq:intro-Sigma-n} also shows that the
limiting field can be analyzed spectrally. The covariance operator is a
periodic convolution operator and is therefore diagonal in the Fourier
basis.

\begin{theorem}[Fourier representation of the limiting covariance]
\label{thm:intro-Fourier-covariance}
Let \(\mathcal K_n\) be the operator on \(L^2([0,1])\) defined by
\begin{equation}
\label{eq:intro-Kn}
(\mathcal K_ng)(x):=a_ng(x)+\int_{\R}b_n(h)g(x+h)\,\dd h,
\end{equation}
where \(g\) is extended \(1\)-periodically to \(\R\). Then
\(\mathcal K_n\) is bounded, positive, and self-adjoint, and
\[
\Sigma_n(f,g)=\langle f,\mathcal K_ng\rangle_{L^2([0,1])}.
\]
The Fourier functions
\[
e_k(x):=e^{2\pi ikx},\qquad k\in\Z,
\]
are eigenfunctions of \(\mathcal K_n\), with eigenvalues
\begin{equation}
\label{eq:intro-theta-nk}
\theta_{n,k}:=a_n+\int_{\R}b_n(h)\cos(2\pi kh)\,\dd h.
\end{equation}
Consequently,
\begin{equation}
\label{eq:intro-Sigma-Fourier}
\Sigma_n(f,f)=\sum_{k\in\Z}\theta_{n,k}|\widehat f(k)|^2.
\end{equation}
\end{theorem}

The Fourier representation has two complementary consequences. First, the
eigenvalues \(\theta_{n,k}\) quantify the strength of the depth-\(n\)
fluctuations at each periodic spatial frequency. The constant mode measures
the total covariance mass, while the nonzero modes describe the response to
oscillatory observation profiles.

Second, the diagonalization yields a prediction principle for the limiting
field. Suppose that the Gaussian coordinates corresponding to finitely many
Fourier modes are observed. The orthogonal projection of \(\xi_n(f)\) onto
these coordinates is its optimal mean-square predictor. The unobserved
remainder is independent of the observed coordinates, and its variance is
the sum of the omitted Fourier contributions
\[
\theta_{n,k}|\widehat f(k)|^2.
\]
Consequently, for large observation windows, integrals of the recursive
field against finitely many Fourier basis functions can be used to estimate
its integral against a more general test function, with an asymptotic
precision determined explicitly by the neglected covariance eigenvalues.

The fixed-depth theorem produces a different limiting Gaussian field for
each \(n\). The next problem is to determine how its covariance grows as the
recursive depth increases. The following assumption guarantees that
sufficiently short chains or branches of Poisson points contribute uniformly
to the recursive field:
\begin{equation}
\label{eq:intro-local-lower-bound}
\varphi(x)\ge c_0,\qquad |x|\le r_0,
\end{equation}
for some \(r_0,c_0>0\).

\begin{theorem}[Growth of the limiting covariance with depth]
\label{thm:intro-covariance-growth}
Assume \eqref{eq:intro-local-lower-bound}. For both recursions,
\begin{equation}
\label{eq:intro-log-growth}
\lim_{n\to\infty}\frac{\log\Sigma_n(1,1)}{d_n\log d_n}=2
\end{equation}
and
\begin{equation}
\label{eq:intro-diagonal-offdiagonal-growth}
\lim_{n\to\infty}\frac{\log a_n}{d_n\log d_n}
=
\lim_{n\to\infty}
\frac{\log\left(\int_{\R}b_n(h)\,\dd h\right)}{d_n\log d_n}
=2.
\end{equation}
If \(f\) satisfies
\[
0<c_f\le f(x)\le C_f<\infty,\qquad x\in\R,
\]
then
\[
\lim_{n\to\infty}\frac{\log\Sigma_n(f,f)}{d_n\log d_n}=2.
\]
\end{theorem}

The diagonal and off-diagonal covariance contributions therefore have the
same leading logarithmic order. The amplification is not generated solely
by repeated self-interaction: covariance between distinct points in the
common Poisson environment remains significant on the scale
\(d_n\log d_n\). The same exponent holds for every positive periodic test
function bounded away from zero.

The two recursion mechanisms differ through the multiplicity \(d_n\). For
the linear recursion, \(d_n=n\), whereas for the nonlinear recursion,
\(d_n=2^n-1\). The covariance therefore obeys the same universal law when
expressed in terms of \(d_n\), but the nonlinear recursion reaches the same
multiplicity at a substantially smaller depth.

Finally, we allow the recursive depth to increase with the observation
window. In this regime, both the multiple-integral order and the dependence
radius change with \(N\). A fixed-depth central limit theorem is therefore
insufficient: the observation window must dominate the increasing recursive
complexity.

\begin{theorem}[Central limit theorem at growing depth]
\label{thm:intro-growing-depth}
Let \(n_N\to\infty\) and set \(D_N:=d_{n_N}\). Suppose that \(f\) is a
bounded \(1\)-periodic function satisfying
\begin{equation}
\label{eq:intro-nonzero-mean}
\int_0^1f(x)\,\dd x\neq0.
\end{equation}
If
\begin{equation}
\label{eq:intro-growing-depth-condition}
D_N\log(D_N+2)=o(\log N),
\end{equation}
then
\begin{equation}
\label{eq:intro-growing-depth-clt}
\frac{Q_{n_N,N}(f)-\E Q_{n_N,N}(f)}
{\sqrt{\Var(Q_{n_N,N}(f))}}
\Longrightarrow\mathcal N(0,1),
\qquad N\to\infty.
\end{equation}
\end{theorem}

This result identifies a regime in which Gaussian approximation remains
valid even though the recursive depth tends to infinity. Because the
multiplicity grows linearly in the linear model and exponentially in the
nonlinear model, the same condition on \(D_N\) leads to different admissible
depth scales for the two recursions.

The nonzero-mean condition makes the constant Fourier mode visible and
provides a uniform lower bound for the variance through the first Poisson
chaos. A more general formulation, expressed through the Fourier response of
the test function, also covers certain mean-zero functions for which a
nonconstant mode remains visible.

Taken together, the results describe three aspects of the recursive system.
The fixed-depth theorem constructs the limiting Gaussian linear family and
identifies its covariance. The Fourier theorem resolves this covariance into
independent spatial modes and provides a basis for predicting unobserved
linear functionals from observed ones. The depth-growth theorem quantifies
the amplification generated by repeated recursion, while the growing-depth
theorem determines when the observation window remains sufficiently large
for Gaussian approximation.

The paper is organized as follows. We first establish locality and moment
properties of the recursive fields. We then identify the limiting covariance
and prove the fixed-depth central limit theorem. The subsequent sections
develop the covariance operator and its Fourier representation, analyze the
growth of the covariance with recursive depth, and prove the growing-depth
central limit theorem.

%===========================================================

%===========================================================
\section{The Gaussian noise associated with a fixed depth}
\label{sec:fixed-depth-Gaussian-noise}
%===========================================================
 
Fix \(n\ge0\). In this section, we construct the Gaussian noise associated
with depth \(n\) from the covariance coefficients \(a_n\) and \(b_n\). We
first show that the corresponding covariance form extends continuously to
\(L^2([0,1])\). We then introduce its covariance operator, derive its Fourier
decomposition, and construct the centered Gaussian linear family
\[
        \{\xi_n(f):f\in L^2([0,1])\}.
\]
The final part of the section formulates the finite-dimensional distributions
of this family and introduces the notions of constant mode, oscillatory
modes, spectral amplification, and spectral nondegeneracy.

The family \(\xi_n\) is naturally interpreted as a generalized Gaussian
field indexed by \(L^2([0,1])\).

%-----------------------------------------------------------
\subsection{The covariance form}
\label{subsec:fixed-depth-covariance-form}
%-----------------------------------------------------------

Let \(\mathcal B_{\mathrm{per}}\) denote the real vector space of bounded
measurable \(1\)-periodic functions on \(\R\). Each
\(f\in\mathcal B_{\mathrm{per}}\) is identified with its restriction to
\([0,1]\).

Fix \(n\ge0\). Recall the coefficients introduced in the Introduction,
\[
        a_n=\lambda\E^0[W_n(0)^2]
\]
and, for \(h\ne0\),
\[
        b_n(h)
        =
        \lambda^2
        \left(
        \E^{0,h}[W_n(0)W_n(h)]
        -
        \bigl(\E^0W_n(0)\bigr)^2
        \right),
        \qquad
        b_n(0):=0.
\]
Their construction from the recursive Poisson fields will be developed in
Section~\ref{sec:covariance-identification}. In particular, we shall prove
there that
\begin{equation}
\label{eq:fixed-depth-proved-covariance-properties}
        a_n\ge0,
        \qquad
        b_n\in L^1(\R),
        \qquad
        b_n(-h)=b_n(h),
\end{equation}
and that the form
\begin{equation}
\label{eq:fixed-depth-Sigma}
        \Sigma_n(f,g)
         :=
        a_n\int_0^1f(x)g(x)\,\dd x
        \quad+
        \int_0^1\int_{\R}
        f(x)g(x+h)b_n(h)\,\dd h\,\dd x
\end{equation}
is nonnegative:
\begin{equation}
\label{eq:fixed-depth-Sigma-nonnegative}
        \Sigma_n(f,f)\ge0
\end{equation}
for every bounded measurable \(1\)-periodic function \(f\).

In the present section, we use these properties to construct and analyze the
Gaussian linear family associated with the fixed depth \(n\).

The coefficient \(a_n\) is the diagonal or local intensity of the
noise and corresponds to pairing a Poisson point with itself. The function
\(b_n(h)\) describes the covariance generated by two distinct weighted
points separated by \(h\). For the recursive models, this dependence is
created by the Poisson input that can influence both depth-\(n\) weights.

Pointwise nonnegativity of \(b_n\), which will hold for the models studied
below, does not by itself imply
\eqref{eq:fixed-depth-Sigma-nonnegative} for sign-changing test functions.
The relevant positivity property is the nonnegativity of the complete
bilinear form.

For \(h\in\R\), let
\begin{equation}
\label{eq:fixed-depth-translation}
        (\tau_hg)(x):=g(x+h),
        \qquad x\in[0,1],
\end{equation}
where the argument is taken modulo \(1\). The operator \(\tau_h\) is unitary
on \(L^2([0,1])\):
\begin{equation}
\label{eq:fixed-depth-translation-unitary}
        \|\tau_hg\|_{L^2([0,1])}
        =
        \|g\|_{L^2([0,1])}.
\end{equation}

\begin{lemma}
%[Bounded extension of the covariance form]
\label{lem:fixed-depth-Sigma-extension}
For \(f,g\in\mathcal B_{\mathrm{per}}\),
\begin{equation}
\label{eq:fixed-depth-Sigma-bound}
        |\Sigma_n(f,g)|
        \le
        \bigl(a_n+\|b_n\|_{L^1(\R)}\bigr)
        \|f\|_{L^2([0,1])}
        \|g\|_{L^2([0,1])}.
\end{equation}
The form \(\Sigma_n\) is symmetric and extends uniquely to a bounded,
symmetric, nonnegative bilinear form on \(L^2([0,1])\).
\end{lemma}

\begin{proof}
Using \eqref{eq:fixed-depth-translation}, the off-diagonal term can be
written as
\[
        \int_{\R}
        b_n(h)
        \langle f,\tau_hg\rangle_{L^2([0,1])}
        \,\dd h.
\]
Therefore, by Cauchy--Schwarz and
\eqref{eq:fixed-depth-translation-unitary},
\begin{align*}
\left|
        \int_{\R}
        b_n(h)
        \langle f,\tau_hg\rangle
        \,\dd h
\right|
&\le
        \int_{\R}
        |b_n(h)|\,
        \|f\|_2\|\tau_hg\|_2
        \,\dd h
\\
&=
        \|b_n\|_1\|f\|_2\|g\|_2.
\end{align*}
The diagonal term satisfies
\[
\left|
        a_n\int_0^1f(x)g(x)\,\dd x
\right|
\le
        a_n\|f\|_2\|g\|_2.
\]
Combining the two estimates gives
\eqref{eq:fixed-depth-Sigma-bound}.

For symmetry, observe that
\[
        \langle f,\tau_hg\rangle
        =
        \langle g,\tau_{-h}f\rangle.
\]
Since \(b_n(-h)=b_n(h)\),
\begin{align*}
\int_{\R}
        b_n(h)\langle f,\tau_hg\rangle\,\dd h
&=
\int_{\R}
        b_n(h)\langle g,\tau_hf\rangle\,\dd h.
\end{align*}
Hence
\[
        \Sigma_n(f,g)=\Sigma_n(g,f).
\]

The space of bounded measurable functions is dense in \(L^2([0,1])\).
Bound \eqref{eq:fixed-depth-Sigma-bound} therefore gives a unique continuous
extension of \(\Sigma_n\) to \(L^2([0,1])\times L^2([0,1])\). If
\(f_j\in\mathcal B_{\mathrm{per}}\) and \(f_j\to f\) in \(L^2([0,1])\),
then continuity and \eqref{eq:fixed-depth-Sigma-nonnegative} imply
\[
        \Sigma_n(f,f)
        =
        \lim_{j\to\infty}\Sigma_n(f_j,f_j)
        \ge0.
\]
Thus the extended form remains nonnegative.
\end{proof}

In the rest of the section, \(\Sigma_n\) denotes this continuous extension to
\(L^2([0,1])\).

%-----------------------------------------------------------
\subsection{The covariance operator and generalized Gaussian element}
\label{subsec:fixed-depth-covariance-operator}
%-----------------------------------------------------------

For \(g\in L^2([0,1])\), define
\begin{equation}
\label{eq:fixed-depth-Tbn}
        (T_{b_n}g)(x)
        :=
        \int_{\R}b_n(h)g(x+h)\,\dd h.
\end{equation}
 By
\eqref{eq:fixed-depth-translation-unitary},
\begin{equation}
\label{eq:fixed-depth-Tbn-bound}
        \|T_{b_n}g\|_2
        \le
        \|b_n\|_1\|g\|_2.
\end{equation}

Define
\begin{equation}
\label{eq:fixed-depth-Kn}
        \mathcal K_n
        :=
        a_nI+T_{b_n},
\end{equation}
or, equivalently,
\begin{equation}
\label{eq:fixed-depth-Kn-pointwise}
        (\mathcal K_ng)(x)
        =
        a_ng(x)
        +
        \int_{\R}b_n(h)g(x+h)\,\dd h.
\end{equation}

\begin{proposition}
%[Covariance operator]
\label{prop:fixed-depth-covariance-operator}
The operator
\[
        \mathcal K_n:
        L^2([0,1])
        \longrightarrow
        L^2([0,1])
\]
is bounded, self-adjoint, and nonnegative. It satisfies
\begin{equation}
\label{eq:fixed-depth-Kn-bound}
        \|\mathcal K_n\|_{\mathrm{op}}
        \le
        a_n+\|b_n\|_{L^1(\R)},
\end{equation}
and
\begin{equation}
\label{eq:fixed-depth-Sigma-operator}
        \Sigma_n(f,g)
        =
        \langle f,\mathcal K_ng\rangle_{L^2([0,1])},
        \qquad
        f,g\in L^2([0,1]).
\end{equation}
\end{proposition}

\begin{proof}
Bound \eqref{eq:fixed-depth-Kn-bound} follows from
\eqref{eq:fixed-depth-Tbn-bound}. For bounded periodic \(f,g\), Fubini's
theorem gives
\[
        \langle f,\mathcal K_ng\rangle
        =
        \Sigma_n(f,g).
\]
Both sides are continuous in \(f\) and \(g\), so
\eqref{eq:fixed-depth-Sigma-operator} extends to all
\(f,g\in L^2([0,1])\).

Since \(\Sigma_n\) is symmetric,
\[
        \langle f,\mathcal K_ng\rangle
        =
        \Sigma_n(f,g)
        =
        \Sigma_n(g,f)
        =
        \langle\mathcal K_nf,g\rangle.
\]
Thus \(\mathcal K_n\) is self-adjoint. Finally,
\[
        \langle f,\mathcal K_nf\rangle
        =
        \Sigma_n(f,f)
        \ge0,
\]
so \(\mathcal K_n\) is nonnegative.
\end{proof}

\begin{example}
%[Action of the covariance operator on a Dirac measure]
\label{ex:fixed-depth-Kn-delta}
Let
\[
        \widetilde b_n(x)
        :=
        \sum_{\ell\in\Z}b_n(x+\ell),
        \qquad x\in\mathbb T,
\]
be the \(1\)-periodization of \(b_n\). The operator \(\mathcal K_n\) has the
 kernel
\begin{equation}
\label{eq:fixed-depth-Cn-kernel}
        C_n(x,y)
        :=
        a_n\delta_x(y)+\widetilde b_n(y-x),
        \qquad x,y\in\mathbb T,
\end{equation}
in the sense that
\begin{equation}
\label{eq:fixed-depth-Kn-kernel-representation}
        (\mathcal K_ng)(x)
        =
        \int_{\mathbb T}C_n(x,y)g(y)\,\dd y,
        \qquad
        g\in C^\infty(\mathbb T).
\end{equation}
Here the first term in \(C_n\) is distributional and represents the
multiplication operator \(a_nI\).

Applying \(\mathcal K_n\) to the Dirac measure at the origin gives
\begin{equation}
\label{eq:fixed-depth-Kn-delta-example}
        \mathcal K_n\delta_0
        =
        C_n(\,\cdot\,,0)
        =
        a_n\delta_0+\widetilde b_n.
\end{equation}
Indeed, for every \(\phi\in C^\infty(\mathbb T)\),
\[
\begin{aligned}
        \langle\mathcal K_n\delta_0,\phi\rangle
        &=
        \langle\delta_0,\mathcal K_n\phi\rangle
=
        (\mathcal K_n\phi)(0)
        \\
        &=
        a_n\phi(0)
        +
        \int_{\mathbb T}\widetilde b_n(y)\phi(y)\,\dd y.
\end{aligned}
\]
Therefore \(\mathcal K_n\delta_0\) consists of a singular part
\(a_n\delta_0\), inherited from the diagonal covariance, and the regular
part \(\widetilde b_n\), generated by the off-diagonal covariance.
\end{example}

The operator \(\mathcal K_n\) represents the covariance of the depth-\(n\)
Gaussian noise through \eqref{eq:fixed-depth-Sigma-operator}. Since it is
bounded, self-adjoint, and nonnegative, the spectral theorem gives a unique
bounded nonnegative square root \(\mathcal K_n^{1/2}\).

We use the standard Hilbert-space construction of a centered Gaussian
linear functional; see Dorogovtsev~\cite[Chapter~1]{DorogovtsevStochasticAnalysisHilbert}.
Let \(G\) be an isonormal Gaussian process over \(L^2([0,1])\), that is,
a centered Gaussian linear map
\[
        G:L^2([0,1])\longrightarrow L^2(\Omega)
\]
satisfying
\[
        \E[G(h)G(q)]
        =
        \langle h,q\rangle_{L^2([0,1])},
        \qquad h,q\in L^2([0,1]).
\]
Define
\begin{equation}
\label{eq:fixed-depth-xi-construction}
        \xi_n(f)
        :=
        G(\mathcal K_n^{1/2}f),
        \qquad f\in L^2([0,1]).
\end{equation}
Then \(\{\xi_n(f):f\in L^2([0,1])\}\) is a generalized Gaussian element
 and
\begin{equation}
\label{eq:fixed-depth-Gaussian-covariance}
        \E[\xi_n(f)\xi_n(g)]
        =
        \left\langle
        \mathcal K_n^{1/2}f,
        \mathcal K_n^{1/2}g
        \right\rangle_{L^2([0,1])}
        =
        \langle f,\mathcal K_ng\rangle_{L^2([0,1])}
        =
        \Sigma_n(f,g).
\end{equation}
Thus \(\Sigma_n\) determines the finite-dimensional distributions of the
centered Gaussian family.

Thus the pair \((a_n,b_n)\) determines the full law of the Gaussian noise at
depth \(n\).

%-----------------------------------------------------------

\subsection{Fourier decomposition and Sobolev realization}
\label{subsec:fixed-depth-Fourier}
%-----------------------------------------------------------

We now diagonalize the covariance operator in the
space \(L^2([0,1])\). This representation will also be used later to compare
the full covariance spectrum with the first-chaos contribution.

For \(k\ge1\), set
\[
        c_k(x):=\sqrt2\cos(2\pi kx),
        \qquad
        s_k(x):=\sqrt2\sin(2\pi kx),
\]
and set \(c_0(x):=1\). Then
\[
        \{c_0\}\cup\{c_k,s_k:k\ge1\}
\]
is an orthonormal basis of \(L^2([0,1])\). For
\(f\in L^2([0,1])\), write
\[
        f_0:=\langle f,c_0\rangle_{L^2},
        \qquad
        f_k^{c}:=\langle f,c_k\rangle_{L^2},
        \qquad
        f_k^{s}:=\langle f,s_k\rangle_{L^2},
        \quad k\ge1.
\]

\begin{proposition}
%[Fourier diagonalization]
\label{prop:fixed-depth-Fourier-diagonalization}
For \(k\ge0\), define
\begin{equation}
\label{eq:fixed-depth-theta-nk}
        \theta_{n,k}
        :=
        a_n+
        \int_{\R}b_n(h)\cos(2\pi kh)\,\dd h .
\end{equation}
Then
\[
        \mathcal K_nc_0=\theta_{n,0}c_0,
\]
and, for every \(k\ge1\),
\[
        \mathcal K_nc_k=\theta_{n,k}c_k,
        \qquad
        \mathcal K_ns_k=\theta_{n,k}s_k .
\]
Moreover,
\begin{equation}
\label{eq:fixed-depth-theta-properties}
        \theta_{n,k}\ge0,
        \qquad
        \lim_{k\to\infty}\theta_{n,k}=a_n,
\end{equation}
and
\begin{equation}
\label{eq:fixed-depth-Kn-exact-norm}
        \|\mathcal K_n\|_{\mathrm{op}}
        =
        \sup_{k\ge0}\theta_{n,k}.
\end{equation}
\end{proposition}

\begin{proof}
We first compute the action of \(T_{b_n}\) on the real Fourier basis. Since
\(b_n\) is even, the sine terms vanish after integration. Hence, for
\(k\ge1\),
\[
\begin{aligned}
        (T_{b_n}c_k)(x)
        &=
        \sqrt2
        \int_{\R}
        b_n(h)\cos(2\pi k(x+h))\,\dd h        \\
        &=
        c_k(x)
        \int_{\R}b_n(h)\cos(2\pi kh)\,\dd h.
\end{aligned}
\]
Therefore \(\mathcal K_nc_k=\theta_{n,k}c_k\). Similarly,
\[
\begin{aligned}
        (T_{b_n}s_k)(x)
        &=
        \sqrt2
        \int_{\R}
        b_n(h)\sin(2\pi k(x+h))\,\dd h        \\
        &=
        s_k(x)
        \int_{\R}b_n(h)\cos(2\pi kh)\,\dd h,
\end{aligned}
\]
and so \(\mathcal K_ns_k=\theta_{n,k}s_k\). Finally, since \(c_0\equiv1\),
\[
        \mathcal K_nc_0
        =
        \left(
        a_n+\int_{\R}b_n(h)\,\dd h
        \right)c_0
        =
        \theta_{n,0}c_0.
\]

The nonnegativity of \(\mathcal K_n\) gives
\[
        \theta_{n,0}
        =
        \langle c_0,\mathcal K_nc_0\rangle_{L^2}
        \ge0,
        \qquad
        \theta_{n,k}
        =
        \langle c_k,\mathcal K_nc_k\rangle_{L^2}
        \ge0,
        \quad k\ge1.
\]
Since \(b_n\in L^1(\R)\), 
\[
        \int_{\R}b_n(h)\cos(2\pi kh)\,\dd h
        \longrightarrow0,
        \qquad k\to\infty.
\]
Thus \(\theta_{n,k}\to a_n\).

It remains to identify the operator norm. Let
\[
        g=g_0c_0+\sum_{k\ge1}(g_k^{c}c_k+g_k^{s}s_k)
\]
be the Fourier expansion of \(g\in L^2([0,1])\). The preceding
diagonalization gives
\[
        \mathcal K_ng
        =
        \theta_{n,0}g_0c_0
        +
        \sum_{k\ge1}
        \theta_{n,k}
        (g_k^{c}c_k+g_k^{s}s_k),
\]
with convergence in \(L^2([0,1])\). Hence
\[
        \|\mathcal K_ng\|_2^2
        =
        \theta_{n,0}^2g_0^2
        +
        \sum_{k\ge1}
        \theta_{n,k}^2
        \bigl((g_k^{c})^2+(g_k^{s})^2\bigr)
        \le
        \left(\sup_{k\ge0}\theta_{n,k}\right)^2
        \|g\|_2^2.
\]
Therefore
\[
        \|\mathcal K_n\|_{\mathrm{op}}
        \le
        \sup_{k\ge0}\theta_{n,k}.
\]
The reverse inequality follows by testing \(\mathcal K_n\) on the normalized
basis vectors \(c_0,c_k,s_k\). This proves
\eqref{eq:fixed-depth-Kn-exact-norm}.
\end{proof}

The diagonalization gives a direct interpretation of the covariance form. If
\[
        f=f_0c_0+\sum_{k\ge1}(f_k^{c}c_k+f_k^{s}s_k),
        \qquad
        g=g_0c_0+\sum_{k\ge1}(g_k^{c}c_k+g_k^{s}s_k),
\]
then
\begin{equation}
\label{eq:fixed-depth-Sigma-Fourier-bilinear}
        \Sigma_n(f,g)
        =
        \theta_{n,0}f_0g_0
        +
        \sum_{k\ge1}
        \theta_{n,k}
        \bigl(
        f_k^{c}g_k^{c}
        +
        f_k^{s}g_k^{s}
        \bigr).
\end{equation}
Thus \(\theta_{n,k}\) is the variance coefficient associated with the
\(k\)-th real periodic frequency. The case \(k=0\) corresponds to the
constant mode, and
\begin{equation}
\label{eq:fixed-depth-constant-mode}
        \theta_{n,0}
        =
        a_n+\int_{\R}b_n(h)\,\dd h
        =
        \Sigma_n(1,1).
\end{equation}
For \(k\ge1\), the cosine and sine components have the same eigenvalue
\(\theta_{n,k}\). Hence each nonzero frequency contributes through the
two-dimensional subspace spanned by \(c_k\) and \(s_k\).

\begin{corollary}
%[Spectral nondegeneracy]
\label{cor:fixed-depth-spectral-nondegeneracy}
For \(f\in L^2([0,1])\), one has
\[
        \Sigma_n(f,f)>0
\]
if and only if the Fourier expansion of \(f\) has a nonzero component in at
least one eigenspace with strictly positive eigenvalue. Equivalently, either
\[
        f_0\ne0
        \quad\text{and}\quad
        \theta_{n,0}>0,
\]
or there exists \(k\ge1\) such that
\[
        (f_k^{c},f_k^{s})\ne(0,0)
        \quad\text{and}\quad
        \theta_{n,k}>0.
\]
\end{corollary}

\begin{proof}
Taking \(g=f\) in \eqref{eq:fixed-depth-Sigma-Fourier-bilinear} gives
\[
        \Sigma_n(f,f)
        =
        \theta_{n,0}f_0^2
        +
        \sum_{k\ge1}
        \theta_{n,k}
        \bigl((f_k^{c})^2+(f_k^{s})^2\bigr).
\]
All terms are nonnegative by
\eqref{eq:fixed-depth-theta-properties}. Therefore the sum is strictly
positive exactly when at least one nonzero Fourier component of \(f\)
belongs to an eigenspace with strictly positive eigenvalue.
\end{proof}

We now identify the regularity of the limiting Gaussian field. Although
\(\mathcal K_n\) is a bounded covariance operator on \(L^2(\mathbb T)\),
this does not imply that the corresponding Gaussian field is
\(L^2(\mathbb T)\)-valued. The diagonal term \(a_nI\) produces a
white-noise component, and the natural realization space is therefore a
negative Sobolev space.

We write 
\[ \mathbb T:=\R/\Z \] 
and denote by \(\mathcal D'(\mathbb T)\) the space of periodic distributions, that is, the topological dual of \(C^\infty(\mathbb T)\).
 
Let \(\Delta_{\mathbb T}\) denote the periodic Laplacian, acting on smooth periodic functions by 
\[ \Delta_{\mathbb T}\phi=\phi''. \]
 Thus, for the real Fourier basis introduced above, 
\[ -\Delta_{\mathbb T}c_0=0, \qquad -\Delta_{\mathbb T}c_k=4\pi^2k^2c_k, \qquad -\Delta_{\mathbb T}s_k=4\pi^2k^2s_k, \quad k\ge1. \] 

we define, for \(r\in\R\), \[ H^r(\mathbb T) := \left\{ u\in\mathcal D'(\mathbb T): \|u\|_{H^r(\mathbb T)}<\infty \right\}, \]
where
\begin{equation} \label{eq:fixed-depth-Hr-norm} \|u\|_{H^r(\mathbb T)}^2 := |u_0|^2 + \sum_{k\ge1} (1+4\pi^2k^2)^r \bigl( |u_k^{c}|^2+|u_k^{s}|^2 \bigr). 
\end{equation}
(see, for example, McLean~\cite[Chapter~3]{McLean2000StronglyElliptic} or Taylor~\cite[Chapter~4]{Taylor2011PDEI} for Sobolev spaces of real order, and Shubin~\cite[Chapter~I]{Shubin2001PDO} for the spectral definition on compact manifolds)

Here
\[
        u_0:=\langle u,c_0\rangle,
        \qquad
        u_k^{c}:=\langle u,c_k\rangle,
        \qquad
        u_k^{s}:=\langle u,s_k\rangle,
        \quad k\ge1.
\]
For \(u\in L^2(\mathbb T)\), these are the usual \(L^2\)-inner products.
For \(u\in\mathcal D'(\mathbb T)\), the same notation denotes the
distributional pairing with the corresponding smooth periodic test
function.

\begin{lemma}
%[Sobolev realization of the fixed-depth Gaussian field]
\label{lem:fixed-depth-Sobolev-realization}
Assume \(a_n>0\). The Gaussian linear family
\[
        \{\xi_n(f):f\in L^2(\mathbb T)\},
        \qquad
        \E[\xi_n(f)\xi_n(g)]
        =
        \Sigma_n(f,g),
\]
does not admit an \(L^2(\mathbb T)\)-valued realization. However, for
every \(s>\frac12\), it admits an \(H^{-s}(\mathbb T)\)-valued centered
Gaussian realization \(\Xi_n\). More precisely, for
\(f,g\in H^s(\mathbb T)\),
\[
        \E\bigl[
        \langle\Xi_n,f\rangle
        \langle\Xi_n,g\rangle
        \bigr]
        =
        \Sigma_n(f,g),
\]
where \(\langle\cdot,\cdot\rangle\) denotes the duality between
\(H^{-s}(\mathbb T)\) and \(H^s(\mathbb T)\). Moreover,
\begin{equation}
\label{eq:fixed-depth-Hminus-s-second-moment}
        \E\|\Xi_n\|_{H^{-s}(\mathbb T)}^2
        =
        \theta_{n,0}
        +
        2\sum_{k\ge1}
        \theta_{n,k}
        (1+4\pi^2k^2)^{-s}
        <\infty .
\end{equation}
The map \(f\mapsto\langle\Xi_n,f\rangle\), initially defined on
\(H^s(\mathbb T)\), extends uniquely in \(L^2(\Omega)\) to all
\(f\in L^2(\mathbb T)\), and this extension is the original Gaussian
family \(f\mapsto\xi_n(f)\).
\end{lemma}

\begin{proof}
By Proposition~\ref{prop:fixed-depth-Fourier-diagonalization},
\[
        \theta_{n,k}
        =
        a_n+\int_{\R}b_n(h)\cos(2\pi kh)\,\dd h .
\]
Since \(b_n\in L^1(\R)\), the Riemann--Lebesgue lemma gives
\[
        \theta_{n,k}\longrightarrow a_n>0.
\]

Suppose first that an \(L^2(\mathbb T)\)-valued centered Gaussian
realization \(X_n\) existed with covariance operator \(\mathcal K_n\).
Using the orthonormal basis
\[
        \{c_0\}\cup\{c_k,s_k:k\ge1\},
\]
we would have, by monotone convergence,
\[
\begin{aligned}
        \E\|X_n\|_{L^2(\mathbb T)}^2
        &=
        \E|\langle X_n,c_0\rangle|^2
        +
        \sum_{k\ge1}
        \E|\langle X_n,c_k\rangle|^2
        +
        \sum_{k\ge1}
        \E|\langle X_n,s_k\rangle|^2        \\
        &=
        \theta_{n,0}
        +
        2\sum_{k\ge1}\theta_{n,k}.
\end{aligned}
\]
The last series diverges because \(\theta_{n,k}\to a_n>0\). Hence no
\(L^2(\mathbb T)\)-valued realization exists.

We now construct the realization in \(H^{-s}(\mathbb T)\). Let
\[
        Z_0,\quad Z_k^{c},\quad Z_k^{s},\qquad k\ge1,
\]
be independent centered real Gaussian random variables with
\[
        \E Z_0^2=\theta_{n,0},
        \qquad
        \E (Z_k^{c})^2=\E (Z_k^{s})^2=\theta_{n,k}.
\]
Define the random Fourier series
\begin{equation}
\label{eq:fixed-depth-Xin-real-series}
        \Xi_n
        :=
        Z_0c_0+
        \sum_{k\ge1}
        \bigl(
        Z_k^{c}c_k+Z_k^{s}s_k
        \bigr).
\end{equation}
For \(s>\frac12\),
\[
\begin{aligned}
        \E\|\Xi_n\|_{H^{-s}(\mathbb T)}^2
        &=
        \theta_{n,0}
        +
        \sum_{k\ge1}
        (1+4\pi^2k^2)^{-s}
        \Bigl(
        \E (Z_k^{c})^2+\E (Z_k^{s})^2
        \Bigr)        \\
        &=
        \theta_{n,0}
        +
        2\sum_{k\ge1}
        \theta_{n,k}
        (1+4\pi^2k^2)^{-s}.
\end{aligned}
\]
The sequence \((\theta_{n,k})_{k\ge0}\) is bounded because
\[
        0\le \theta_{n,k}\le \|\mathcal K_n\|_{\mathrm{op}}.
\]
Therefore
\[
        \sum_{k\ge1}
        \theta_{n,k}
        (1+4\pi^2k^2)^{-s}
        \le
        \|\mathcal K_n\|_{\mathrm{op}}
        \sum_{k\ge1}
        (1+4\pi^2k^2)^{-s}
        <\infty
\]
whenever \(s>\frac12\). Thus
\eqref{eq:fixed-depth-Xin-real-series} converges in
\(L^2(\Omega;H^{-s}(\mathbb T))\).

For a real trigonometric polynomial
\[
        f=f_0c_0+\sum_{k\ge1}(f_k^{c}c_k+f_k^{s}s_k),
\]
we have
\[
        \langle\Xi_n,f\rangle
        =
        Z_0f_0+
        \sum_{k\ge1}
        \bigl(
        Z_k^{c}f_k^{c}
        +
        Z_k^{s}f_k^{s}
        \bigr).
\]
Hence, for trigonometric polynomials \(f,g\),
\[
\begin{aligned}
        \E[
        \langle\Xi_n,f\rangle
        \langle\Xi_n,g\rangle
        ]
        &=
        \theta_{n,0}f_0g_0
        +
        \sum_{k\ge1}
        \theta_{n,k}
        \bigl(
        f_k^{c}g_k^{c}
        +
        f_k^{s}g_k^{s}
        \bigr)        \\
        &=
        \Sigma_n(f,g),
\end{aligned}
\]
where the last equality is
\eqref{eq:fixed-depth-Sigma-Fourier-bilinear}. By density, this identity
extends to all \(f,g\in H^s(\mathbb T)\).

Finally, since \(\Sigma_n\) is bounded on
\(L^2(\mathbb T)\times L^2(\mathbb T)\), the Gaussian linear map
defined on trigonometric polynomials extends uniquely in \(L^2(\Omega)\)
to all \(f\in L^2(\mathbb T)\). Its covariance is still \(\Sigma_n\), so
this extension is exactly the Gaussian family \(f\mapsto\xi_n(f)\).
\end{proof}

Thus the fixed-depth limit is not a pointwise Gaussian process and not an
\(L^2\)-valued random function when \(a_n>0\). It is, however, a genuine
Hilbert-space-valued Gaussian random distribution in \(H^{-s}(\mathbb T)\)
for every \(s>\frac12\). The threshold is the same as for one-dimensional
white noise; the off-diagonal Palm term modifies the covariance but does not
remove the white-noise component produced by the diagonal term \(a_nI\).

\begin{example}[Depth zero: the white-noise benchmark]
\label{ex:fixed-depth-zero}
At depth zero, \(W_0\equiv1\). Hence
\[
 a_0=\lambda,
 \qquad
 b_0(h)=0,
 \qquad
 \mathcal K_0=\lambda I,
\]
and therefore
\[
 \theta_{0,k}=\lambda,
 \qquad k\in\Z.
\]
Consequently, the limiting Gaussian family at depth zero is periodic white
noise with intensity \(\lambda\): for \(f,g\in L^2(\mathbb T)\),
\[
 \E[\xi_0(f)\xi_0(g)]
 =\lambda\int_0^1 f(x)g(x)\,\dd x.
\]
Since the eigenvalues do not decay as \(|k|\to\infty\), this example also
illustrates why the limiting field has no \(L^2(\mathbb T)\)-valued
realization, whereas it admits an \(H^{-s}(\mathbb T)\)-valued realization
for every \(s>\frac12\).
\end{example}
%-----------------------------------------------------------
\subsection{Finite-dimensional distributions}
\label{subsec:fixed-depth-fdd}
%-----------------------------------------------------------

For every \(m\ge1\) and every
\(f_1,\ldots,f_m\in L^2([0,1])\), the vector
\begin{equation}
\label{eq:fixed-depth-Gaussian-vector}
        \bigl(
        \xi_n(f_1),\ldots,\xi_n(f_m)
        \bigr)
\end{equation}
is centered Gaussian with covariance matrix
\begin{equation}
\label{eq:fixed-depth-Gaussian-matrix}
        \left(
        \Sigma_n(f_i,f_j)
        \right)_{1\le i,j\le m}.
\end{equation}
The nonnegativity of \(\Sigma_n\) guarantees that
\eqref{eq:fixed-depth-Gaussian-matrix} is nonnegative definite.

The Gaussian family is defined on all of \(L^2([0,1])\). The probabilistic
central limit theorem proved later is formulated first on the dense class
\(\mathcal B_{\mathrm{per}}\) of bounded measurable periodic test functions,
for which the Poisson statistics \(Q_{n,N}(f)\) have already been defined.

For \(f\in\mathcal B_{\mathrm{per}}\), set
\begin{equation}
\label{eq:fixed-depth-random-functional}
        \mathcal Z_{n,N}(f)
        :=
        \frac{
        Q_{n,N}(f)-\E Q_{n,N}(f)
        }{
        \sqrt{|\Lambda_N|}
        }.
\end{equation}
The fixed-depth central limit theorem will establish that, for every
\(m\ge1\) and every
\(f_1,\ldots,f_m\in\mathcal B_{\mathrm{per}}\),
\begin{equation}
\label{eq:fixed-depth-fdd-convergence}
        \bigl(
        \mathcal Z_{n,N}(f_1),\ldots,
        \mathcal Z_{n,N}(f_m)
        \bigr)
        \Longrightarrow
        \bigl(
        \xi_n(f_1),\ldots,
        \xi_n(f_m)
        \bigr).
\end{equation}

\begin{remark}[Approximation from finitely many Fourier coordinates]
\label{rem:fixed-depth-Fourier-prediction}

Let \(f\in L^2([0,1])\) have the Fourier expansion
\[
f=f_0c_0+\sum_{j\geq1}\bigl(f_j^{c}c_j+f_j^{s}s_j\bigr),
\]
and let
\[
P_Kf:=f_0c_0+\sum_{j=1}^{K}\bigl(f_j^{c}c_j+f_j^{s}s_j\bigr).
\]
If the integrals of the weighted point measure \(W_n(v)\mu(\dd v)\) against \(c_0,c_1,s_1,\ldots,c_K,s_K\) are known, then linearity gives
\[
Q_{n,N}(P_Kf)=f_0Q_{n,N}(c_0)+\sum_{j=1}^{K}\bigl(f_j^{c}Q_{n,N}(c_j)+f_j^{s}Q_{n,N}(s_j)\bigr).
\]

For the limiting Gaussian functional,
\[
\xi_n(f)=\xi_n(P_Kf)+\xi_n(f-P_Kf).
\]
By the diagonalization of the covariance operator, \(\xi_n(P_Kf)\) is measurable with respect to
\[
\mathcal G_{n,K}:=\sigma\bigl(\xi_n(c_0),\xi_n(c_j),\xi_n(s_j):1\leq j\leq K\bigr),
\]
whereas \(\xi_n(f-P_Kf)\) is independent of \(\mathcal G_{n,K}\). Hence
\[
\E\bigl[\xi_n(f)\mid\mathcal G_{n,K}\bigr]=\xi_n(P_Kf).
\]
Thus \(\xi_n(P_Kf)\) is the optimal mean-square predictor of \(\xi_n(f)\) from the first \(K\) Fourier coordinates. Its prediction error is
\begin{equation}
\label{eq:fixed-depth-Fourier-prediction-error}
\E\bigl[|\xi_n(f)-\xi_n(P_Kf)|^2\bigr]=\sum_{j>K}\theta_{n,j}\bigl((f_j^{c})^2+(f_j^{s})^2\bigr).
\end{equation}
In particular,
\[
\E\bigl[|\xi_n(f)-\xi_n(P_Kf)|^2\bigr]\leq \|\mathcal K_n\|_{\mathrm{op}}\|f-P_Kf\|_{L^2([0,1])}^2,
\]
and the right-hand side converges to zero as \(K\to\infty\).

Writing \(\overline Q_{n,N}(g):=Q_{n,N}(g)-\E[Q_{n,N}(g)]\), the finite-dimensional central limit theorem therefore suggests, for large \(N\), the approximation
\[
Q_{n,N}(f)\approx \E[Q_{n,N}(f)]+f_0\overline Q_{n,N}(c_0)+\sum_{j=1}^{K}\bigl(f_j^{c}\overline Q_{n,N}(c_j)+f_j^{s}\overline Q_{n,N}(s_j)\bigr).
\]
For each fixed \(K\), the normalized approximation error satisfies
\[
\frac{\overline Q_{n,N}(f)-\overline Q_{n,N}(P_Kf)}{\sqrt{|\Lambda_N|}}\Longrightarrow \mathcal N\!\left(0,\sum_{j>K}\theta_{n,j}\bigl((f_j^{c})^2+(f_j^{s})^2\bigr)\right).
\]
Therefore, finitely many Fourier integrals provide an asymptotic estimate of the integral against an arbitrary test function, while \eqref{eq:fixed-depth-Fourier-prediction-error} quantifies the loss of precision caused by truncation.
\end{remark}

\section{Proof of the fixed-depth Gaussian limit}
\label{sec:recursive-models}
\label{sec:covariance-identification}
\label{sec:fixed-depth-clt}

This section proves Theorem~\ref{thm:intro-fixed-depth}.  The two recursive
models, the unified field notation \(W_n\), the multiplicity \(d_n\), and the
statistics \(Q_{n,N}(f)\) were introduced in
Section~\ref{sec:introduction}.  We therefore do not redefine the models.
Instead, we establish the three ingredients required for the theorem: the
locality and moment properties of the recursive fields, the Palm
representation of the limiting covariance, and the central limit theorem for
the resulting finitely dependent block sequence.  Throughout the section,
\(n\geq0\) is fixed and \(W_n\) denotes either \(F_n\) or \(A_n\).

\subsection{Local structure of the recursive fields}
\label{subsec:fixed-depth-local-structure}

Let \(\mathcal N_{\mathrm{lf}}(\R)\) be the space of locally finite counting
measures on \(\R\), equipped with the evaluation sigma-field.  For
\(\zeta\in\mathcal N_{\mathrm{lf}}(\R)\), let \(F_n^\zeta\) and
\(A_n^\zeta\) denote the fields obtained from
\eqref{eq:intro-linear-recursion} and
\eqref{eq:intro-nonlinear-recursion}, respectively, after replacing the
Poisson measure \(\mu\) by \(\zeta\).  We write
\[
        W_n^\zeta=F_n^\zeta
        \qquad\text{or}\qquad
        W_n^\zeta=A_n^\zeta,
\]
and \(W_n=W_n^\mu\).  This notation will be used only when the dependence on
the underlying configuration must be explicit.  Recall also that
\[
        d_n=n
        \quad\text{for the linear recursion},
        \qquad
        d_n=2^n-1
        \quad\text{for the nonlinear recursion}.
\]
Set
\begin{equation}
\label{eq:kernel-uniform-bound}
        M:=\|\varphi\|_\infty.
\end{equation}

\begin{lemma}[Pathwise local-count bound]
\label{lem:pathwise-local-count}
For either recursion, every \(\zeta\in\mathcal N_{\mathrm{lf}}(\R)\),
\(v\in\R\), and \(n\geq0\),
\begin{equation}
\label{eq:pathwise-weight-bound}
        0\leq W_n^\zeta(v)
        \leq
        M^{d_n}\,
        \zeta([v-nR,v+nR])^{d_n},
\end{equation}
with the convention \(x^0=1\).
\end{lemma}

\begin{proof}
The assertion is immediate for \(n=0\).  Suppose it holds at depth \(n\),
and put
\[
        Z:=\zeta([v-(n+1)R,v+(n+1)R]).
\]
Whenever \(|u-v|\leq R\),
\[
        [u-nR,u+nR]
        \subset
        [v-(n+1)R,v+(n+1)R].
\]
For the linear recursion, the induction hypothesis therefore gives
\[
\begin{aligned}
        F_{n+1}^\zeta(v)
        &\leq
        \int_{[v-R,v+R]}
        M\,M^{d_n}Z^{d_n}\,\zeta(\dd u)  \\
        &\leq M^{d_n+1}Z^{d_n+1}.
\end{aligned}
\]
Since \(d_{n+1}=d_n+1\) in the linear model, this is
\eqref{eq:pathwise-weight-bound} at depth \(n+1\).

For the nonlinear recursion, the induction hypothesis yields
\[
        A_n^\zeta(v)\leq M^{d_n}Z^{d_n}
\]
and
\[
        \int_{[v-R,v+R]}
        \varphi(u-v)A_n^\zeta(u)\,\zeta(\dd u)
        \leq M^{d_n+1}Z^{d_n+1}.
\]
Multiplication gives
\[
        A_{n+1}^\zeta(v)
        \leq M^{2d_n+1}Z^{2d_n+1}.
\]
Here \(d_{n+1}=2d_n+1\), which completes the induction.
\end{proof}

\begin{proposition}[Measurability, locality, and moments]
\label{prop:recursive-fields-measurable}
For every fixed \(n\), the maps
\[
        (v,\zeta)\longmapsto F_n^\zeta(v),
        \qquad
        (v,\zeta)\longmapsto A_n^\zeta(v)
\]
are jointly measurable and finite on
\(\R\times\mathcal N_{\mathrm{lf}}(\R)\).  Moreover,
\begin{equation}
\label{eq:recursive-field-locality}
        W_n^\zeta(v)
        \text{ depends only on }
        \zeta\big|_{[v-nR,v+nR]}.
\end{equation}
Under the Poisson law, the restrictions of these fields to bounded sets have
moments of every finite order.  The same conclusion remains valid after the
addition of finitely many deterministic atoms to \(\mu\).
\end{proposition}

\begin{proof}
Joint measurability is proved by induction.  It is immediate at depth zero.
Assume that \((u,\zeta)\mapsto W_n^\zeta(u)\) is jointly measurable.  Then
\[
        (v,u,\zeta)
        \longmapsto
        \varphi(u-v)W_n^\zeta(u)\1_{\{|u-v|\leq R\}}
\]
is jointly measurable.  Integration against the counting-measure coordinate
preserves measurability; see
\cite[Chapter~2]{LastPenrose2017}.  In the nonlinear recursion,
multiplication by \(A_n^\zeta(v)\) preserves this property.  Compact support
of \(\varphi\) and local finiteness of \(\zeta\) imply finiteness of every
recursive sum.

The locality assertion follows from the same induction.  At depth zero there
is nothing to prove.  At depth \(n+1\), only atoms
\(u\in[v-R,v+R]\) enter the convolution.  For each such \(u\), the induction
hypothesis restricts \(W_n^\zeta(u)\) to
\([u-nR,u+nR]\), which is contained in
\([v-(n+1)R,v+(n+1)R]\).  In the nonlinear model, the additional factor
\(A_n^\zeta(v)\) is measurable with respect to the smaller interval
\([v-nR,v+nR]\).

Finally, let \(K\subset\R\) be bounded.  All intervals
\([v-nR,v+nR]\), \(v\in K\), are contained in the bounded enlargement
\(K^{+nR}\).  If \(q\) deterministic atoms are added to \(\mu\),
Lemma~\ref{lem:pathwise-local-count} gives
\[
        \sup_{v\in K}W_n^{\mu+\sum_{j=1}^q\delta_{x_j}}(v)
        \leq
        M^{d_n}\bigl(\mu(K^{+nR})+q\bigr)^{d_n}.
\]
Every finite moment of the Poisson count on the right-hand side is finite.
\end{proof}

For \(t\in\R\), define the translated configuration by
\[
        (\vartheta_t\zeta)(B):=\zeta(B-t).
\]
The recursive definitions and a change of variables give
\begin{equation}
\label{eq:translation-covariance-recursive-fields}
        W_n^{\vartheta_t\zeta}(v+t)=W_n^\zeta(v),
        \qquad t,v\in\R.
\end{equation}
Thus a stationary input produces a stationary recursive field.

We shall also use the natural order on counting measures:
\(\zeta\leq\zeta'\) if
\(\zeta(B)\leq\zeta'(B)\) for every Borel set \(B\).

\begin{lemma}[Monotonicity]
\label{lem:recursive-field-monotonicity}
If \(\zeta\leq\zeta'\), then
\[
        W_n^\zeta(v)\leq W_n^{\zeta'}(v),
        \qquad v\in\R,
\]
for both recursions and every \(n\geq0\).
\end{lemma}

\begin{proof}
The proof is by induction.  In the linear recursion, both the nonnegative
integrand and the integrating counting measure increase.  In the nonlinear
recursion, the same argument applies to the convolution factor, while the
additional factor \(A_n(v)\) also increases by the induction hypothesis.
\end{proof}

\subsection{The limiting covariance}
\label{subsec:fixed-depth-covariance-identification}

Recall the coefficients \(a_n\) and \(b_n\) from
\eqref{eq:intro-an}--\eqref{eq:intro-bn}.  For the Poisson process,
Slivnyak's theorem identifies the one-point and two-point Palm expectations
with deterministic insertions:
\begin{equation}
\label{eq:slivnyak-insertions}
        \E^x H(\mu)=\E H(\mu+\delta_x),
        \qquad
        \E^{x,y}H(\mu)
        =\E H(\mu+\delta_x+\delta_y),
        \quad x\neq y.
\end{equation}
This representation is used below only to evaluate local functionals of the
configuration.

\begin{lemma}[Positivity of the diagonal coefficient]
\label{lem:diagonal-coefficient-positive}
For both recursive models,
\begin{equation}
\label{eq:an-self-loop-lower}
        a_n\geq\lambda\alpha_0^{2d_n}>0.
\end{equation}
\end{lemma}

\begin{proof}
Under the one-point Palm distribution at the origin, the configuration
contains an atom at \(0\).  In the linear recursion, retaining only this atom
at every convolution step yields
\(F_n(0)\geq\alpha_0^n\).  In the nonlinear recursion,
\[
        A_{n+1}(0)\geq \alpha_0 A_n(0)^2.
\]
Since \(A_0(0)=1\) and \(d_{n+1}=2d_n+1\), induction gives
\(A_n(0)\geq\alpha_0^{d_n}\).  Squaring and taking the one-point Palm
expectation proves the claim.
\end{proof}

\begin{lemma}[Support, symmetry, and integrability of the covariance density]
\label{lem:fixed-depth-offdiag-kernel-range}
The function \(b_n\) is measurable and even.  It satisfies
\begin{equation}
\label{eq:bn-support}
        b_n(h)=0,
        \qquad |h|>2nR,
\end{equation}
and is bounded on its support.  Consequently,
\(b_n\in L^1(\R)\).
\end{lemma}

\begin{proof}
For \(h\neq0\), the two-point Palm distribution is represented by the
configuration \(\mu+\delta_0+\delta_h\).  If \(|h|>2nR\), the dependence
intervals
\[
        [-nR,nR]
        \qquad\text{and}\qquad
        [h-nR,h+nR]
\]
are disjoint.  Moreover, neither inserted atom belongs to the dependence
interval of the field evaluated at the other location.  By
Proposition~\ref{prop:recursive-fields-measurable}, the two field values are
therefore functions of independent restrictions of \(\mu\), together with
their respective deterministic atoms.  Hence
\[
\begin{aligned}
        \E^{0,h}[W_n(0)W_n(h)]
        &=
        \E[W_n^{\mu+\delta_0}(0)]
        \E[W_n^{\mu+\delta_h}(h)] \\
        &=\bigl(\E^0W_n(0)\bigr)^2,
\end{aligned}
\]
where translation covariance was used in the second identity.  This proves
\eqref{eq:bn-support}.

Measurability follows from joint measurability in
Proposition~\ref{prop:recursive-fields-measurable}.  To prove evenness,
translate the configuration in the two-point Palm expectation by \(-h\)
and interchange the two factors.  This gives
\[
        \E^{0,h}[W_n(0)W_n(h)]
        =
        \E^{0,-h}[W_n(0)W_n(-h)].
\]

It remains to establish boundedness on the support.  The case \(n=0\) is
immediate.  For \(|h|\leq2nR\), both relevant dependence intervals are
contained in \([-3nR-1,3nR+1]\).  With
\[
        Z_n:=\mu([-3nR-1,3nR+1]),
\]
Lemma~\ref{lem:pathwise-local-count} gives, under the two-point Palm
distribution,
\[
        W_n(0)\vee W_n(h)
        \leq M^{d_n}(Z_n+2)^{d_n}.
\]
The required two-point moments are therefore bounded uniformly over
\(|h|\leq2nR\), because every moment of \(Z_n\) is finite.  The deterministic
term \((\E^0W_n(0))^2\) is finite by
Proposition~\ref{prop:recursive-fields-measurable}.  Thus \(b_n\) is bounded
and compactly supported.
\end{proof}

\begin{proposition}[Nonnegativity of the covariance density]
\label{prop:bn-nonnegative-association}
For every \(n\geq0\) and \(h\in\R\),
\[
        b_n(h)\geq0.
\]
\end{proposition}

\begin{proof}
The assertion at \(h=0\) follows from the convention \(b_n(0)=0\).  Fix
\(h\neq0\) and define, as functionals of the background Poisson process,
\[
        X_h:=W_n^{\mu+\delta_0+\delta_h}(0),
        \qquad
        Y_h:=W_n^{\mu+\delta_0+\delta_h}(h).
\]
By Lemma~\ref{lem:recursive-field-monotonicity}, both functionals are
increasing.  The Poisson Harris--FKG inequality
\cite[Theorem~20.4]{LastPenrose2017}, applied first to bounded truncations and
then followed by monotone convergence, yields
\[
        \E[X_hY_h]\geq \E[X_h]\,\E[Y_h].
\]
Monotonicity and translation covariance further imply
\[
        \E[X_h]\geq \E W_n^{\mu+\delta_0}(0)=\E^0W_n(0),
        \qquad
        \E[Y_h]\geq \E W_n^{\mu+\delta_h}(h)=\E^0W_n(0).
\]
Using \eqref{eq:slivnyak-insertions}, we obtain
\[
        \E^{0,h}[W_n(0)W_n(h)]
        \geq
        \bigl(\E^0W_n(0)\bigr)^2,
\]
which is equivalent to \(b_n(h)\geq0\).
\end{proof}

\begin{proposition}[Finite-window covariance and its limit]
\label{prop:fixed-depth-covariance-kernel}
For bounded \(1\)-periodic functions \(f,g\),
\begin{align}
\Cov(Q_{n,N}(f),Q_{n,N}(g))
&=a_n\int_{\Lambda_N}f(v)g(v)\,\dd v \notag\\
&\quad+
\int_{\Lambda_N}\int_{\Lambda_N}
        f(v)g(u)b_n(u-v)\,\dd u\,\dd v.
\label{eq:fixed-depth-covariance-kernel-finiteN}
\end{align}
Consequently,
\begin{equation}
\label{eq:fixed-depth-covariance-limit}
        \frac{1}{2N}
        \Cov(Q_{n,N}(f),Q_{n,N}(g))
        \longrightarrow
        \Sigma_n(f,g),
\end{equation}
where \(\Sigma_n\) is the covariance form defined in
\eqref{eq:fixed-depth-Sigma}.
\end{proposition}

\begin{proof}
The one-point Mecke formula gives
\begin{equation}
\label{eq:mean-QnN-fixed-depth}
        \E Q_{n,N}(f)
        =
        \lambda\E^0W_n(0)
        \int_{\Lambda_N}f(v)\,\dd v.
\end{equation}
Indeed, apply the formula to
\(H(v,\mu)=f(v)W_n^\mu(v)\1_{\Lambda_N}(v)\), and then use translation
covariance.

For the second moment, decompose the ordered pairs of atoms into the diagonal
part and the factorial, off-diagonal part.  The one-point and two-point Mecke
formulas \cite[Chapter~4]{LastPenrose2017} yield
\[
\begin{aligned}
        \E[Q_{n,N}(f)Q_{n,N}(g)]
        &=
        \lambda\int_{\Lambda_N}
        f(v)g(v)\E^v[W_n(v)^2] \,\dd v \\
        &\quad+
        \lambda^2\int_{\Lambda_N^2}
        f(v)g(u)\E^{v,u}[W_n(v)W_n(u)]
        \,\dd u\,\dd v.
\end{aligned}
\]
The first term equals
\(a_n\int_{\Lambda_N}fg\).  For \(u\neq v\), translation covariance gives
\[
        \E^{v,u}[W_n(v)W_n(u)]
        =
        \E^{0,u-v}[W_n(0)W_n(u-v)].
\]
Subtracting the product of the means in
\eqref{eq:mean-QnN-fixed-depth} and using the definition of \(b_n\) proves
\eqref{eq:fixed-depth-covariance-kernel-finiteN}.  The diagonal
\(\{u=v\}\) is Lebesgue-null in the off-diagonal integrals.

For the limit, the diagonal contribution divided by \(2N\) equals
\[
        a_n\int_0^1f(x)g(x)\,\dd x
\]
exactly.  In the second term, set \(h=u-v\).  Its normalized form is
\[
        \int_\R b_n(h)
        \frac{1}{2N}
        \int_{\Lambda_N\cap(\Lambda_N-h)}
        f(v)g(v+h)\,\dd v\,\dd h.
\]
For fixed \(h\), the inner integrand is bounded and \(1\)-periodic in
\(v\), and therefore
\[
        \frac{1}{2N}
        \int_{\Lambda_N\cap(\Lambda_N-h)}
        f(v)g(v+h)\,\dd v
        \longrightarrow
        \int_0^1f(x)g(x+h)\,\dd x.
\]
The absolute value is bounded by
\(\|f\|_\infty\|g\|_\infty\).  Since \(b_n\in L^1(\R)\), dominated
convergence proves \eqref{eq:fixed-depth-covariance-limit}.  In particular,
\(\Sigma_n(f,f)\geq0\), because it is the limit of normalized variances.
\end{proof}

The following consequence will be used when the depth tends to infinity.

\begin{corollary}[Comparison with the constant mode]
\label{cor:comparison-with-constant-mode}
If
\[
        0\leq m_f\leq f(x)\leq M_f<\infty,
        \qquad x\in\R,
\]
then
\[
        m_f^2\Sigma_n(1,1)
        \leq
        \Sigma_n(f,f)
        \leq
        M_f^2\Sigma_n(1,1).
\]
\end{corollary}

\begin{proof}
Apply the pointwise bounds to the two terms of
\eqref{eq:fixed-depth-Sigma}, using \(a_n\geq0\) and \(b_n\geq0\), and
then use \eqref{eq:fixed-depth-constant-mode}.
\end{proof}

\subsection{Block decomposition and proof of Theorem~\ref{thm:intro-fixed-depth}}
\label{subsec:fixed-depth-block-clt}

Let
\[
        I_k:=[k,k+1),
        \qquad k\in\Z,
\]
and define
\begin{equation}
\label{eq:block-contribution-fixed-n}
\begin{aligned}
        B_{n,k}(f)
        &:=
        \int_{I_k}f(v)W_n(v)\,\mu(\dd v),\\
        X_{n,k}(f)
        &:=
        B_{n,k}(f)-\E B_{n,k}(f).
\end{aligned}
\end{equation}
The next estimate is quantitative in the depth and will also be used in the
growing-depth analysis.

\begin{proposition}[Uniform block moment estimate]
\label{prop:block-moment-estimate}
For every \(p\geq1\), there exists
\[
        C_{p,f}=C(p,\lambda,R,M,\|f\|_\infty)<\infty
\]
such that
\begin{equation}
\label{eq:block-moment-estimate}
        \|X_{n,k}(f)\|_p
        \leq
        C_{p,f}^{\,d_n+1}(1+d_n)^{d_n+1},
        \qquad n\geq0,
        \quad k\in\Z.
\end{equation}
\end{proposition}

The proof is given in the section on moment bounds and technical estimates.
For the fixed-depth argument below, only the consequence
\(X_{n,0}(f)\in L^2\) is needed.

Set
\begin{equation}
\label{eq:m-n-fixed-depth}
        m(n):=\lceil2Rn\rceil+2.
\end{equation}

\begin{proposition}[Stationary finite-range block decomposition]
\label{thm:block-decomposition-finite-dependence}
For fixed \(n\), the sequence
\((X_{n,k}(f))_{k\in\Z}\) is strictly stationary and
\(m(n)\)-dependent.  Moreover,
\begin{equation}
\label{eq:centered-Q-block-decomposition-fixed-n}
        Q_{n,N}(f)-\E Q_{n,N}(f)
        =
        \sum_{k=-N}^{N-1}X_{n,k}(f).
\end{equation}
\end{proposition}

\begin{proof}
By Proposition~\ref{prop:recursive-fields-measurable}, the block variable
\(B_{n,k}(f)\) is measurable with respect to the restriction of \(\mu\) to
\[
        J_{n,k}:=[k-nR,k+1+nR].
\]
If two finite sets of block indices are separated by more than \(m(n)\), the
corresponding unions of the intervals \(J_{n,k}\) are disjoint.  Independent
increments of the Poisson process therefore imply independence of the two
families of block variables.  This is \(m(n)\)-dependence.

For an integer \(r\), translation covariance
\eqref{eq:translation-covariance-recursive-fields}, stationarity of \(\mu\),
and the identity \(f(v+r)=f(v)\) show that
\[
        (B_{n,k_1+r}(f),\ldots,B_{n,k_q+r}(f))
        \stackrel{d}{=}
        (B_{n,k_1}(f),\ldots,B_{n,k_q}(f)).
\]
Thus the block sequence is strictly stationary.

To obtain the exact sum, use the half-open representative
\[
        \Lambda_N^\circ:=[-N,N)
        =\bigcup_{k=-N}^{N-1}I_k.
\]
The union is disjoint.  Since the intensity measure is diffuse,
\(\mu(\{N\})=0\) almost surely, and replacing \(\Lambda_N\) by
\(\Lambda_N^\circ\) does not change \(Q_{n,N}(f)\).  Hence
\[
        Q_{n,N}(f)=\sum_{k=-N}^{N-1}B_{n,k}(f)
        \qquad\text{almost surely}.
\]
Subtracting expectations proves
\eqref{eq:centered-Q-block-decomposition-fixed-n}.
\end{proof}

For later use, if \(n=n_N\), define
\begin{equation}
\label{eq:mN-definition}
        D_N:=d_{n_N},
        \qquad
        m_N:=m(n_N),
        \qquad
        X_{N,k}(f):=X_{n_N,k}(f).
\end{equation}
Then each row is stationary and \(m_N\)-dependent, and
\begin{equation}
\label{eq:growing-depth-block-sum}
        Q_{n_N,N}(f)-\E Q_{n_N,N}(f)
        =
        \sum_{k=-N}^{N-1}X_{N,k}(f).
\end{equation}

\begin{proof}[Proof of Theorem~\ref{thm:intro-fixed-depth}]
Lemmas~\ref{lem:fixed-depth-offdiag-kernel-range} and
Proposition~\ref{prop:bn-nonnegative-association} show that \(b_n\) is even,
nonnegative, and supported on \([-2nR,2nR]\).  Proposition
\ref{prop:fixed-depth-covariance-kernel} identifies the limiting covariance
with \(\Sigma_n\).  Section~\ref{sec:fixed-depth-Gaussian-noise} then shows
that \(\Sigma_n\) extends to a bounded symmetric nonnegative form on
\(L^2([0,1])\) and constructs the centered Gaussian family \(\xi_n\) with
this covariance.  It remains to prove the finite-dimensional convergence.

Fix \(q\geq1\), functions
\(f_1,\ldots,f_q\in\mathcal B_{\mathrm{per}}\), and coefficients
\(t_1,\ldots,t_q\in\R\).  Put
\[
        f_t:=\sum_{j=1}^q t_jf_j
\]
and
\[
        S_N
        :=Q_{n,N}(f_t)-\E Q_{n,N}(f_t)
        =\sum_{k=-N}^{N-1}X_{n,k}(f_t),
\]
where the last identity follows from
Proposition~\ref{thm:block-decomposition-finite-dependence}.  Let
\(s_N^2:=\Var(S_N)\).  By
\eqref{eq:fixed-depth-covariance-limit},
\begin{equation}
\label{eq:fixed-depth-variance-per-length}
        \frac{s_N^2}{2N}
        \longrightarrow
        \Sigma_n(f_t,f_t).
\end{equation}

Assume first that \(\Sigma_n(f_t,f_t)>0\).  Reindex the row by
\[
        Y_{N,j}:=X_{n,j-N-1}(f_t),
        \qquad 1\leq j\leq2N.
\]
The row is \(m(n)\)-dependent and centered.  We apply the
\(m\)-dependent Lindeberg criterion of
\cite[Theorem~1.1]{JansonPratelliRigo2024}.  For every
\(\varepsilon>0\), its Lindeberg quantity is, by stationarity,
\[
        \frac{2Nm(n)}{s_N^2}
        \E\!\left[
        X_{n,0}(f_t)^2
        \1_{\{|X_{n,0}(f_t)|>
        \varepsilon s_N/m(n)\}}
        \right].
\]
By \eqref{eq:fixed-depth-variance-per-length},
\[
        \frac{2Nm(n)}{s_N^2}
        \longrightarrow
        \frac{m(n)}{\Sigma_n(f_t,f_t)},
        \qquad
        \frac{\varepsilon s_N}{m(n)}\longrightarrow\infty.
\]
The random variable \(X_{n,0}(f_t)^2\) is integrable by
Proposition~\ref{prop:block-moment-estimate}.  Dominated convergence therefore
verifies the Lindeberg condition, and
\[
        \frac{S_N}{s_N}
        \Longrightarrow
        \mathcal N(0,1).
\]
Together with \eqref{eq:fixed-depth-variance-per-length} and Slutsky's
theorem, this yields
\[
        \frac{S_N}{\sqrt{2N}}
        \Longrightarrow
        \mathcal N\bigl(0,\Sigma_n(f_t,f_t)\bigr).
\]

If \(\Sigma_n(f_t,f_t)=0\), then
\[
        \E\left|\frac{S_N}{\sqrt{2N}}\right|^2
        =\frac{s_N^2}{2N}\longrightarrow0,
\]
so the same conclusion holds with the degenerate centered Gaussian law.
Thus every linear combination of
\[
        \left(
        \frac{Q_{n,N}(f_j)-\E Q_{n,N}(f_j)}{\sqrt{2N}}
        \right)_{j=1}^q
\]
converges to the corresponding centered Gaussian linear combination with
variance \(\Sigma_n(f_t,f_t)\).  The Cram\'er--Wold device and
\eqref{eq:fixed-depth-Gaussian-covariance} prove
\eqref{eq:intro-fixed-depth-fdd}.
\end{proof}

% \begin{remark}
% The convergence proved above is finite-dimensional convergence of the random
% linear functionals.  Convergence in an \(H^{-s}(\mathbb T)\) topology would
% require a separate tightness argument.
%\end{remark}
\section{Covariance growth and Gaussian limits at increasing depth}
\label{sec:depth-evolution}
\label{sec:growing-depth-clt}

Sections~\ref{sec:fixed-depth-Gaussian-noise} and
\ref{sec:fixed-depth-clt} identify the limiting Gaussian family for each
fixed recursion depth.  We now study how its covariance changes with the
depth and then allow the depth to increase with the observation window.
These questions are linked: the increasing-depth central limit theorem
requires a variance lower bound that remains effective while the dependence
range grows.  We obtain such a bound from the first Poisson chaos of the
statistic.

This section proves Theorems~\ref{thm:intro-covariance-growth} and
\ref{thm:intro-growing-depth}, which were stated in the Introduction.  We do
not repeat their statements.  We first establish the covariance-growth
estimates, then construct an add-one response operator and derive spectral
variance bounds, and finally verify a triangular-array central limit theorem.
All notation from Sections~\ref{sec:introduction}--\ref{sec:fixed-depth-clt}
remains in force.

\subsection{Covariance growth with the recursion depth}
\label{subsec:log-growth-constant-mode-covariance}

\begin{assumption}[Local positivity]
\label{ass:local-positivity-kernel}
The kernel satisfies the lower bound
\eqref{eq:intro-local-lower-bound} introduced in the statement of
Theorem~\ref{thm:intro-covariance-growth}.
\end{assumption}

For a merely measurable kernel, this condition is stronger than
\(\varphi(0)>0\).  For example, \(\varphi=\1_{\{0\}}\) satisfies
\(\varphi(0)=1\) but is not bounded away from zero on any neighborhood of
the origin.

\begin{lemma}
\label{lem:poisson-high-moment-log}
If \(Z\) is Poisson with fixed mean \(\beta>0\), then, for integers
\(q\to\infty\),
\[
 \log\E Z^q=q\log q+o(q\log q).
\]
\end{lemma}

\begin{proof}
Since \(x^q\le q!e^x\) for \(x\ge0\),
\[
 \log\E Z^q\le\log(q!)+\beta(e-1)=q\log q+O(q).
\]
For the reverse bound, let \(k_q=\lfloor q/\log q\rfloor\). Then
\[
 \E Z^q\ge k_q^q e^{-\beta}\frac{\beta^{k_q}}{k_q!}.
\]
Using \(\log(k_q!)\le k_q\log k_q\) gives
\[
 \log\E Z^q\ge q\log q-q\log\log q-O(q),
\]
which proves the claim.
\end{proof}

We first proof the asymptotic of covariance for constant function $f=1$. As stated in Theorem~\ref{thm:intro-covariance-growth},
\begin{equation*}
        \lim_{n\to\infty}
        \frac{
        \log\Sigma_n(1,1)
        }{
        d_n\log d_n
        }
        =2,
\end{equation*}
\begin{proof}
By Proposition~\ref{prop:bn-nonnegative-association},
\(\Sigma_n(1,1)\ge a_n\). Under the one-point Palm law, write
\(\zeta=\mu+\delta_0\), let
\[
 I=[-r_0/2,r_0/2],\qquad Z=\mu(I),
\]
and note that \(\zeta(I)=Z+1\). We claim, pathwise, that
\begin{equation}
\label{eq:local-cloud-weight-lower}
 W_n^\zeta(v)\ge c_0^{d_n}\zeta(I)^{d_n},
 \qquad v\in I.
\end{equation}
This avoids any measurability issue associated with an infimum over \(v\).
The claim is immediate at depth zero. If it holds at depth \(n\), then for
every \(v\in I\) and every \(u\in I\), \(|u-v|\le r_0\), so
\[
\begin{split}
 F_{n+1}^\zeta(v)
 &\ge\int_I c_0\,c_0^{d_n}\zeta(I)^{d_n}\,\zeta(\dd u)\\
 &=c_0^{d_n+1}\zeta(I)^{d_n+1}.
\end{split}
\]
For the nonlinear recursion,
\[
\begin{split}
 A_{n+1}^\zeta(v)
 &\ge
 \bigl(c_0^{d_n}\zeta(I)^{d_n}\bigr)
 \int_Ic_0\bigl(c_0^{d_n}\zeta(I)^{d_n}\bigr)\,\zeta(\dd u)\\
 &=c_0^{2d_n+1}\zeta(I)^{2d_n+1}.
\end{split}
\]
The identities \(d_{n+1}=d_n+1\) and \(d_{n+1}=2d_n+1\) prove
\eqref{eq:local-cloud-weight-lower} by induction. In particular,
\[
 W_n^{\mu+\delta_0}(0)
 \ge c_0^{d_n}(Z+1)^{d_n}
 \ge c_0^{d_n}Z^{d_n}.
\]
Therefore
\[
 a_n\ge\lambda c_0^{2d_n}\E Z^{2d_n}.
\]
Lemma~\ref{lem:poisson-high-moment-log} yields
\[
 \liminf_{n\to\infty}
 \frac{\log\Sigma_n(1,1)}{d_n\log d_n}\ge2.
\]

For the upper bound, stationarity and \(m(n)\)-dependence give
\[
\begin{split}
 \Var\!\left(\sum_{k=-N}^{N-1}X_{n,k}(1)\right)
 &=\sum_{k,\ell=-N}^{N-1}
   \Cov(X_{n,k}(1),X_{n,\ell}(1))\\
 &=\sum_{\substack{k,\ell=-N\\|k-\ell|\le m(n)}}^{N-1}
   \Cov(X_{n,k}(1),X_{n,\ell}(1))\\
 &\le2N(2m(n)+1)\|X_{n,0}(1)\|_2^2.
\end{split}
\]
The last step uses Cauchy--Schwarz and stationarity; pairs outside the
displayed band have zero covariance by independence.
Divide by \(2N\), let \(N\to\infty\), and apply Proposition
\ref{prop:block-moment-estimate} with \(p=2\). Since
\(m(n)\le C(1+d_n)\),
\[
 \Sigma_n(1,1)
 \le C^{2d_n+3}(1+d_n)^{2d_n+3}.
\]
Thus the corresponding limsup is at most \(2\), proving
\eqref{eq:intro-log-growth}.
\end{proof}

\subsection{Add-one response, spectral lower bounds, and off-diagonal growth}
\label{subsec:first-chaos-response-operator}

The constant-mode estimate alone is insufficient for increasing depth,
because a test function may have little or no constant Fourier component.
We therefore identify the first-chaos response of the statistic to the
insertion of a single Poisson point.  The resulting translation-invariant
operator yields a mode-by-mode lower bound for the variance.  We then return
to the local positivity assumption and complete the proof of
Theorem~\ref{thm:intro-covariance-growth} by estimating the off-diagonal
covariance mass.

For a square-integrable Poisson functional \(F=F(\mu)\), define the add-one difference
\[
 D_zF(\mu):=F(\mu+\delta_z)-F(\mu).
\]
For the finite-window statistic, put
\begin{equation}
\label{eq:first-chaos-kernel-finiteN}
 \Delta_{n,N}^f(z):=\E[D_zQ_{n,N}(f)].
\end{equation}
Write \(W_n^{(z)}:=W_n^{\mu+\delta_z}\).  The corresponding change
of the weighted point measure is the random signed measure
\begin{equation}
\label{eq:Gamma-nz-definition}
 \Gamma_{n,z}(\dd v):=
 W_n^{(z)}(v)(\mu+\delta_z)(\dd v)-W_n(v)\mu(\dd v).
\end{equation}
Then
\begin{equation}
\label{eq:DzQ-Gamma-representation}
 D_zQ_{n,N}(f)=\int_{\Lambda_N}f(v)\Gamma_{n,z}(\dd v).
\end{equation}

\begin{lemma}
\label{lem:Gamma-local-support}
The measure \(\Gamma_{n,z}\) is supported on
\([z-nR,z+nR]\).
\end{lemma}

\begin{proof}
If \(|v-z|>nR\), the inserted point does not belong to the dependence
interval of \(W_n(v)\), so \(W_n^{(z)}(v)=W_n(v)\). The only new atom is
at \(z\), which lies in the stated interval.
\end{proof}

Define the finite signed measure
\begin{equation}
\label{eq:gamma-n-response-measure}
 \gamma_n(A):=\E[\Gamma_{n,0}(A)].
\end{equation}
It is supported on \([-nR,nR]\). Finiteness follows from
Lemma~\ref{lem:pathwise-local-count}. Indeed, from
\[
\begin{split}
 \Gamma_{n,0}(\dd v)
 &=\bigl(W_n^{(0)}(v)-W_n(v)\bigr)\mu(\dd v)
   +W_n^{(0)}(0)\delta_0(\dd v)
\end{split}
\]
and the triangle inequality,
\[
\begin{split}
 |\Gamma_{n,0}|(\R)
 &\le\int_{[-nR,nR]}
    \bigl(W_n^{(0)}(v)+W_n(v)\bigr)\,\mu(\dd v)
    +W_n^{(0)}(0).
\end{split}
\]
Every dependence interval appearing on the right is contained in
\([-2nR,2nR]\). Put \(Z_n:=\mu([-2nR,2nR])+1\). The local-count bound is
at most a deterministic constant times a fixed power of \(Z_n\). Since
\(Z_n\) has moments of every finite order,
\(\E|\Gamma_{n,0}|(\R)<\infty\), and \(\gamma_n\) is a finite signed
measure.

For \(A\in\mathcal B(\R)\), translation covariance of the inserted
configuration gives
\[
 \E[\Gamma_{n,z}(z+A)]=\E[\Gamma_{n,0}(A)]=\gamma_n(A).
\]
Integrating a bounded measurable function against this identity yields
\begin{equation}
\label{eq:expected-response-gamma}
 \E\!\left[\int_\R f(v)\Gamma_{n,z}(\dd v)\right]
 =\int_\R f(z+s)\gamma_n(\dd s).
\end{equation}

\begin{definition}
\label{def:first-chaos-response-operator}
For bounded periodic \(f\), define
\begin{equation}
\label{eq:Pn-interior-response}
 (\mathcal P_nf)(z):=\int_\R f(z+s)\gamma_n(\dd s).
\end{equation}
\end{definition}

The operator extends to \(L^2(\T)\), with
\begin{equation}
\label{eq:Pn-L2-bound}
 \|\mathcal P_nf\|_2\le|\gamma_n|(\R)\|f\|_2.
\end{equation}
Indeed, for bounded \(f\),
\[
 |(\mathcal P_nf)(z)|
 \le\int_\R|f(z+s)|\,|\gamma_n|(\dd s),
\]
and Minkowski's integral inequality gives
\eqref{eq:Pn-L2-bound}. Density then gives the extension. It commutes with
translations because
\[
 \mathcal P_n(\tau_tf)(z)
 =\int_\R f(z+t+s)\gamma_n(\dd s)
 =(\tau_t\mathcal P_nf)(z).
\]
It is therefore diagonal in the complex Fourier basis
\(e_k(x)=e^{2\pi ikx}\):
\begin{equation}
\label{eq:pnk-gamma-representation}
 \mathcal P_ne_k=p_{n,k}e_k,\qquad
 p_{n,k}:=\int_\R e^{2\pi iks}\gamma_n(\dd s).
\end{equation}
Indeed,
\[
 (\mathcal P_ne_k)(z)
 =e^{2\pi ikz}\int_\R e^{2\pi iks}\gamma_n(\dd s).
\]
If
\(\widehat f(k)=\int_0^1f(x)e^{-2\pi ikx}\,\dd x\), Parseval's identity
gives
\begin{equation}
\label{eq:Pn-Fourier-expansion}
 \|\mathcal P_nf\|_2^2
 =\sum_{k\in\Z}|p_{n,k}|^2|\widehat f(k)|^2.
\end{equation}

\begin{example}[Add-one response at depths zero and one]
\label{ex:first-depth-response}
At depth zero, \(W_0\equiv1\), so insertion of an atom at \(z\) gives
\[
 \Gamma_{0,z}=\delta_z,
 \qquad
 \gamma_0=\delta_0,
 \qquad
 p_{0,k}=1,
 \quad k\in\Z.
\]
Thus the first-chaos lower bound is an equality at depth zero:
\[
 \theta_{0,k}=\lambda|p_{0,k}|^2=\lambda.
\]

The linear and nonlinear recursions also coincide at depth one:
\[
 W_1(v)=\int_\R\varphi(u-v)\,\mu(\dd u).
\]
After inserting an atom at \(z\),
\[
 W_1^{(z)}(v)-W_1(v)=\varphi(z-v).
\]
Using the definition of \(\Gamma_{n,z}\) in
\eqref{eq:Gamma-nz-definition}, we obtain
\[
\begin{aligned}
 \Gamma_{1,z}(\dd v)
 &=\varphi(z-v)\,\mu(\dd v) \\
 &\quad+
 \left(
   \int_\R\varphi(u-z)\,\mu(\dd u)+\alpha_0
 \right)\delta_z(\dd v).
\end{aligned}
\]
Setting \(z=0\) and applying the Campbell formula gives
\begin{equation*}
 \gamma_1(\dd s)
 =\lambda\varphi(-s)\,\dd s
 +\left(
    \lambda\int_\R\varphi(u)\,\dd u+\alpha_0
  \right)\delta_0(\dd s).
\end{equation*}
Define
\[
 \widehat\varphi(k)
 :=\int_\R\varphi(u)e^{-2\pi iku}\,\dd u,
 \qquad k\in\Z.
\]
Then
\begin{equation*}
 p_{1,k}
 =\alpha_0
  +\lambda\int_\R\varphi(u)\,\dd u
  +\lambda\widehat\varphi(k),
 \qquad k\in\Z.
\end{equation*}
In particular,
\[
 p_{1,0}
 =\alpha_0+2\lambda\int_\R\varphi(u)\,\dd u>0.
\]
If \(\varphi\) is not symmetric, the multipliers \(p_{1,k}\) may be
complex for \(k\ne0\); the spectral lower bound depends on their modulus:
\[
 \theta_{1,k}\ge\lambda|p_{1,k}|^2.
\]
\end{example}

\begin{lemma}[Interior response identity]
\label{lem:Delta-equals-Pn-on-interior}
If \(z\in[-N+nR,N-nR)\), then
\begin{equation}
\label{eq:Delta-equals-Pn-on-interior}
 \Delta_{n,N}^f(z)=(\mathcal P_nf)(z).
\end{equation}
\end{lemma}

\begin{proof}
For such \(z\), the support of \(\Gamma_{n,z}\) lies in \(\Lambda_N\).
Equations \eqref{eq:DzQ-Gamma-representation} and
\eqref{eq:expected-response-gamma} then give the identity.
\end{proof}

\medskip
\noindent\textbf{First-chaos lower bound.}

For \(F\in L^2\) and \(q\ge1\), let
\[
 T_qF(z_1,\ldots,z_q):=\E[D^q_{z_1,\ldots,z_q}F].
\]
The Poisson Fock-space identity states
\begin{equation}
\label{eq:Fock-variance}
 \Var(F)=\sum_{q=1}^{\infty}\frac1{q!}
 \|T_qF\|_{L^2((\lambda\,\dd z)^q)}^2;
\end{equation}
see \cite[Theorem~18.6]{LastPenrose2017}.

\begin{proposition}[First-chaos lower bounds]
\label{prop:first-chaos-spectral-lower-bound}
For bounded periodic \(f\),
\begin{equation}
\label{eq:first-chaos-lower-bound-finiteN}
 \Var(Q_{n,N}(f))\ge
 \lambda\int_\R|\Delta_{n,N}^f(z)|^2\,\dd z,
\end{equation}
and
\begin{equation}
\label{eq:first-chaos-spectral-lower-bound}
 \Sigma_n(f,f)\ge\lambda\|\mathcal P_nf\|_2^2
 =\lambda\sum_{k\in\Z}|p_{n,k}|^2|\widehat f(k)|^2.
\end{equation}
In particular,
\begin{equation}
\label{eq:theta-first-chaos-comparison}
 \theta_{n,k}\ge\lambda|p_{n,k}|^2,\qquad k\in\Z.
\end{equation}
\end{proposition}

\begin{proof}
Proposition~\ref{prop:block-moment-estimate} shows that \(Q_{n,N}(f)\in L^2\), so
\eqref{eq:Fock-variance} applies with \(F=Q_{n,N}(f)\). Every chaos
contribution is nonnegative, and its first term is
\[
 \|T_1F\|_{L^2(\lambda\,\dd z)}^2
 =\lambda\int_\R|\E D_zQ_{n,N}(f)|^2\,\dd z
 =\lambda\int_\R|\Delta_{n,N}^f(z)|^2\,\dd z.
\]
This proves \eqref{eq:first-chaos-lower-bound-finiteN}.

Restrict the last integral to the interior interval and use
Lemma~\ref{lem:Delta-equals-Pn-on-interior}:
\[
 \frac{\Var(Q_{n,N}(f))}{2N}
 \ge\frac{\lambda}{2N}
 \int_{-N+nR}^{N-nR}|(\mathcal P_nf)(z)|^2\,\dd z.
\]
The integrand is \(1\)-periodic. The left-hand side converges to
\(\Sigma_n(f,f)\), by Proposition
\ref{prop:fixed-depth-covariance-kernel}, while the right-hand side
converges to
\(\lambda\int_0^1|(\mathcal P_nf)(x)|^2\,\dd x\). Hence
\[
 \Sigma_n(f,f)\ge\lambda\|\mathcal P_nf\|_2^2.
\]
Substitution of \eqref{eq:Pn-Fourier-expansion} gives the asserted series.

For \(k=0\), apply the last inequality to \(c_0\). For \(k\ge1\), apply
it to the normalized real mode \(c_k\). Since \(\gamma_n\) is real,
\(p_{n,-k}=\overline{p_{n,k}}\), and hence
\[
 \mathcal P_nc_k
 =\frac{p_{n,k}e_k+p_{n,-k}e_{-k}}{\sqrt2},
\qquad
 \|\mathcal P_nc_k\|_2^2=|p_{n,k}|^2.
\]
On the other hand,
\(\Sigma_n(c_k,c_k)=\theta_{n,k}\). This proves
\eqref{eq:theta-first-chaos-comparison}; negative \(k\) follow by symmetry.
\end{proof}

For later use, define
\begin{equation}
\label{eq:rho-n-definition}
 \rho_n(f):=\|\mathcal P_nf\|_2,\qquad
 \rho_n(f)^2=\sum_{k\in\Z}|p_{n,k}|^2|\widehat f(k)|^2.
\end{equation}

\medskip
\noindent\textbf{Positivity and the zero Fourier mode.}

\begin{lemma}
\label{lem:positive-response-measure}
The measure \(\gamma_n\) is positive and
\begin{equation}
\label{eq:pn0-positive}
 p_{n,0}=\gamma_n(\R)\ge\alpha_0^{d_n}>0.
\end{equation}
Consequently, \(\mathcal P_n\) is positivity preserving.
\end{lemma}

\begin{proof}
Lemma~\ref{lem:recursive-field-monotonicity} shows that
\(W_n^{(z)}(v)\ge W_n(v)\). Thus the change of mass at every original
atom is nonnegative, and the inserted atom contributes
\(W_n^{(z)}(z)\delta_z\). Hence \(\Gamma_{n,z}\), and therefore
\(\gamma_n\), is positive.

The inserted atom at zero gives
\(\Gamma_{n,0}(\R)\ge W_n^{(0)}(0)\). Repeatedly retaining only the
self-interaction term gives
\(W_n^{(0)}(0)\ge\alpha_0^{d_n}\) in both recursions; the two
recursion-specific inductions are written in the proof of
Lemma~\ref{lem:diagonal-coefficient-positive}.
Taking expectations proves \eqref{eq:pn0-positive}.
\end{proof}

\begin{corollary}
\label{cor:zero-mode-first-chaos-nondegeneracy}
If \(f\in L^2(\T)\) is real and \(\widehat f(0)\ne0\), then
\begin{equation}
\label{eq:rho-zero-mode-lower-bound}
 \rho_n(f)\ge\alpha_0^{d_n}|\widehat f(0)|,
\end{equation}
and
\begin{equation}
\label{eq:Sigma-zero-mode-lower-bound}
 \Sigma_n(f,f)\ge
 \lambda\alpha_0^{2d_n}|\widehat f(0)|^2.
\end{equation}
\end{corollary}

\begin{proof}
Retain the \(k=0\) term in \eqref{eq:rho-n-definition}, use
\eqref{eq:pn0-positive}, and then apply
\eqref{eq:first-chaos-spectral-lower-bound}.
\end{proof}

Positivity of \(\gamma_n\) does not force \(p_{n,k}\ne0\) for
\(k\ne0\); Fourier transforms of positive finite measures may vanish at
nonzero frequencies. Mean-zero test functions therefore require explicit
control of at least one nonzero response multiplier.

\medskip
\noindent\textbf{A shared local configuration and the off-diagonal covariance.}

To estimate \(\int_\R b_n(h)\,\dd h\), we isolate a local Poisson
configuration that contributes simultaneously to the fields rooted at two
nearby locations.  The following representation makes this common
contribution explicit.

\begin{lemma}[Positive product-measure representation]
\label{lem:positive-multiple-integral-representation}
For every \(n\ge0\) and \(v\in\R\), there is a nonnegative measurable
kernel \(K_{n,v}:\R^{d_n}\to[0,\infty)\) such that
\begin{equation}
\label{eq:positive-multiple-integral-representation}
 W_n^\zeta(v)=\int_{\R^{d_n}}K_{n,v}(\mathbf x)
 \,\zeta^{\otimes d_n}(\dd\mathbf x)
\end{equation}
for every locally finite \(\zeta\). The kernel is a product of \(d_n\)
values of \(\varphi\), one for every edge of the corresponding rooted
integration tree.
\end{lemma}

\begin{proof}
For the linear recursion this is the expansion \eqref{eq:intro-linear-walk-expansion}, with
\[
 K_{n,v}(x_1,\ldots,x_n)
 =\prod_{j=1}^n\varphi(x_j-x_{j-1}),\qquad x_0=v.
\]
Suppose the representation is available for \(A_n\). Substituting it twice
into
\[
 A_{n+1}^\zeta(v)
 =A_n^\zeta(v)\int_\R
   \varphi(u-v)A_n^\zeta(u)\,\zeta(\dd u)
\]
and applying Tonelli's theorem gives one product-measure integral over
\(2d_n+1\) variables. Its kernel is
\[
 K_{n,v}(\mathbf x)\,
 \varphi(u-v)\,
 K_{n,u}(\mathbf y),
\]
which is nonnegative and is a product of \(2d_n+1=d_{n+1}\) kernel
values. This proves the assertion by induction.
\end{proof}

% \begin{lemma}[Product-measure monotonicity]
% \label{lem:product-measure-monotonicity}
% Let \(q\ge1\), let \(H:\R^q\to[0,\infty)\) be measurable, and let some
% coordinates be integrated against fixed positive measures while the remaining
% coordinates are integrated against a locally finite counting measure
% \(\zeta\). Whenever the integral is finite, the resulting functional is
% increasing in \(\zeta\).
% \end{lemma}

% \begin{proof}
% If \(\zeta'\ge\zeta\), write \(\zeta'=\zeta+\nu\) with a positive
% locally finite measure \(\nu\). Expanding every varying factor in the
% product measure gives
% \[
%  (\zeta+\nu)^{\otimes q}
%  =\sum_{\varepsilon\in\{0,1\}^q}
%    \bigotimes_{i=1}^q\zeta_{\varepsilon_i},
%  \qquad \zeta_0=\zeta,\quad\zeta_1=\nu.
% \]
% After fixed coordinate measures are inserted, the same expansion applies
% only to the varying coordinates. The all-zero term is the original integral
% and every other term is nonnegative.
% \end{proof}

Under Assumption~\ref{ass:local-positivity-kernel}, put
\begin{equation}
\label{eq:common-cloud-interval}
 I=[-r_0/2,r_0/2],\qquad
 Z:=\mu(I)\sim\operatorname{Poisson}(\lambda r_0).
\end{equation}

\begin{proposition}
%[Shared-configuration lower bound]
\label{prop:common-cloud-bn-lower-bound}
For either recursion, \(n\ge1\), and \(0<|h|\le r_0/2\),
\begin{equation}
\label{eq:common-cloud-pointwise-lower-bound}
 b_n(h)\ge\lambda^2c_0^{2d_n}\Var(Z^{d_n}).
\end{equation}
Consequently,
\begin{equation}
\label{eq:common-cloud-integral-lower-bound}
 \int_\R b_n(h)\,\dd h
 \ge r_0\lambda^2c_0^{2d_n}\Var(Z^{d_n}).
\end{equation}
\end{proposition}

\begin{proof}
Fix \(h\) in the stated range, write \(d=d_n\), and set
\[
 X=W_n^{\mu+\delta_0+\delta_h}(0),\qquad
 Y=W_n^{\mu+\delta_0+\delta_h}(h),\qquad
 \overline w_n=\E W_n^{\mu+\delta_0}(0).
\]
Then
\[
 \frac{b_n(h)}{\lambda^2}
 =\Cov(X,Y)+\E X\,\E Y-\overline w_n^2
 \ge\Cov(X,Y),
\]
because insertion monotonicity and stationarity give
\(\E X,\E Y\ge\overline w_n\).

Let \(K_v=K_{n,v}\) for \(v=0,h\). If every coordinate belongs to \(I\),
then every edge of either integration tree has length at most \(r_0\);
hence \(K_v\ge c_0^d\) on \(I^d\). To make the required decomposition
explicit, put
\[
 \nu_*=\mu,\qquad \nu_0=\delta_0,\qquad \nu_h=\delta_h.
\]
Expanding the product measure gives
\[
 W_n^{\mu+\delta_0+\delta_h}(v)
 =\sum_{\varepsilon\in\{*,0,h\}^d}
 \int_{\R^d}K_v(\mathbf x)
 \bigotimes_{i=1}^d\nu_{\varepsilon_i}(\dd x_i).
\]
Separate the all-background term and define
\[
 G:=c_0^d\int_{I^d}\mu^{\otimes d}(\dd\mathbf x)
   =c_0^dZ^d
\]
and
\[
\begin{split}
 R_v&:=
 \int_{I^d}(K_v(\mathbf x)-c_0^d)
       \,\mu^{\otimes d}(\dd\mathbf x)\\
 &\quad+\int_{\R^d\setminus I^d}K_v(\mathbf x)
       \,\mu^{\otimes d}(\dd\mathbf x)\\
 &\quad+
 \sum_{\varepsilon\ne(*,\ldots,*)}
 \int_{\R^d}K_v(\mathbf x)
 \bigotimes_{i=1}^d\nu_{\varepsilon_i}(\dd x_i).
\end{split}
\]
Thus \(X=G+R_0\) and \(Y=G+R_h\). Every kernel in the displayed
% remainders is nonnegative. Lemma
% \ref{lem:product-measure-monotonicity}, applied also to terms containing
% fixed measures \(\delta_0,\delta_h\), therefore shows that
\(G,R_0,R_h\) are nonnegative increasing functionals of \(\mu\). They
are square-integrable by
Lemma~\ref{lem:pathwise-local-count}. The Harris--FKG inequality for Poisson processes
now gives
\[
\begin{split}
 \Cov(X,Y)
 &=\Var(G)+\Cov(G,R_h)+\Cov(R_0,G)+\Cov(R_0,R_h)\\
 &\ge\Var(G)=c_0^{2d}\Var(Z^d).
\end{split}
\]
This proves \eqref{eq:common-cloud-pointwise-lower-bound}. Integration over
\(\{h:0<|h|\le r_0/2\}\), whose Lebesgue measure is \(r_0\), proves
\eqref{eq:common-cloud-integral-lower-bound}.
\end{proof}

\begin{lemma}
\label{lem:poisson-power-variance-log}
If \(Z\) is Poisson with fixed positive mean, then
\begin{equation}
\label{eq:poisson-power-variance-log}
 \log\Var(Z^q)=2q\log q+o(q\log q),
 \qquad q\to\infty,\quad q\in\N.
\end{equation}
\end{lemma}

\begin{proof}
The upper bound follows from
\(\Var(Z^q)\le\E Z^{2q}\) and
Lemma~\ref{lem:poisson-high-moment-log}:
\[
 \log\Var(Z^q)\le2q\log(2q)+o(q\log q)
 =2q\log q+o(q\log q).
\]
For the reverse inequality, let \(Z'\) be an independent copy and put
\(k_q=\lfloor q/\log q\rfloor\). The symmetrization identity gives
\[
 \Var(Z^q)=\tfrac12\E(Z^q-(Z')^q)^2
 \ge\tfrac12 k_q^{2q}\Pr\{Z=k_q\}\Pr\{Z'=0\}.
\]
If the Poisson mean is \(\beta>0\), then
\[
 \Pr\{Z=k_q\}\Pr\{Z'=0\}
 =e^{-2\beta}\frac{\beta^{k_q}}{k_q!}.
\]
Therefore, using \(\log(k_q!)\le k_q\log k_q\),
\[
 \log\Var(Z^q)
 \ge2q\log k_q-k_q\log k_q+k_q\log\beta-O(1).
\]
Now
\(\log k_q=\log q-\log\log q+o(1)\) and
\(k_q\log k_q=O(q)\), so
\[
 \log\Var(Z^q)
 \ge2q\log q-2q\log\log q-O(q).
\]
The last two error terms are \(o(q\log q)\), matching the upper bound.
\end{proof}

\begin{proof}[Completion of the proof of
Theorem~\ref{thm:intro-covariance-growth}]
The proof of \eqref{eq:intro-log-growth} gives
\(a_n\ge\lambda c_0^{2d_n}\E Z_0^{2d_n}\) for a fixed-mean Poisson
variable \(Z_0\). Lemma~\ref{lem:poisson-high-moment-log} gives the lower
limit, while
\(a_n\le\Sigma_n(1,1)\) and \eqref{eq:intro-log-growth} give the upper limit.
For the off-diagonal term, Proposition
\ref{prop:common-cloud-bn-lower-bound} and Lemma
\ref{lem:poisson-power-variance-log} give the lower logarithmic bound, while
\[
 0\le\int_\R b_n(h)\,\dd h\le\Sigma_n(1,1)
\]
and \eqref{eq:intro-log-growth} give the upper bound.  Hence
\eqref{eq:intro-diagonal-offdiagonal-growth} holds.

Finally, if \(0<c_f\le f\le C_f<\infty\), Corollary
\ref{cor:comparison-with-constant-mode} gives
\[
 c_f^2\Sigma_n(1,1)
 \le \Sigma_n(f,f)
 \le C_f^2\Sigma_n(1,1).
\]
Taking logarithms and using \eqref{eq:intro-log-growth} proves the remaining
assertion of Theorem~\ref{thm:intro-covariance-growth}.
\end{proof}

The argument shows that the diagonal coefficient and the integrated
off-diagonal covariance have the same leading logarithmic exponent, but it
does not determine their asymptotic ratio.  Since \(d_n=n\) in the linear
model and \(d_n=2^n-1\) in the nonlinear model,
\[
 \log\Sigma_n(1,1)=(2+o(1))n\log n
\]
in the linear case, whereas
\[
 \log\Sigma_n(1,1)=(2\log2+o(1))n2^n
\]
in the nonlinear case.  These are logarithmic equivalents, not ratio
asymptotics for \(\Sigma_n(1,1)\).

\subsection{Gaussian limits when the depth increases with the window}
\label{subsec:triangular-array-criterion}

Let \(n_N\to\infty\).  Recall from \eqref{eq:mN-definition} that
\[
 D_N=d_{n_N},\qquad
 m_N=\lceil2Rn_N\rceil+2,
 \qquad X_{N,k}(f)=X_{n_N,k}(f).
\]
Define
\begin{equation}
\label{eq:growing-depth-variance-definition}
 s_N^2(f):=\Var\bigl(Q_{n_N,N}(f)\bigr),
 \qquad s_N(f)\ge0.
\end{equation}
By \eqref{eq:growing-depth-block-sum}, the centered statistic is a sum of a
stationary \(m_N\)-dependent row.  We first record the triangular-array
criterion used below.

The next result is Utev's triangular-array central limit theorem in the
formulation reproduced by \cite[Theorem~1]{JansonPratelliRigo2024}.  We
use only its qualitative Lindeberg implication.

\begin{theorem}[Utev's criterion]
\label{thm:dependent-type-criterion}
For each \(N\), let \(X_{N,1},\ldots,X_{N,L_N}\) be centered,
square-integrable random variables, let
\[
 S_N=\sum_{j=1}^{L_N}X_{N,j},\qquad s_N^2=\Var(S_N)>0,
\]
and suppose the row is \(m_N\)-dependent, with
\(1\le m_N\le L_N\). If, for every \(\varepsilon>0\),
\begin{equation}
\label{eq:dependent-lindeberg-condition}
 \frac{m_N}{s_N^2}\sum_{j=1}^{L_N}
 \E\!\left[X_{N,j}^2
 \1_{\{|X_{N,j}|>\varepsilon s_N/m_N\}}\right]\longrightarrow0,
\end{equation}
then \(S_N/s_N\Longrightarrow\mathcal N(0,1)\).
\end{theorem}

Here \(m_N\)-dependence has the source's meaning: for each \(j\), the
families with indices at most \(j\) and at least \(j+m_N+1\) are
independent. Our row \((X_{N,k})_{-N\le k<N}\) is put in this form by the
deterministic reindexing \(j=k+N+1\); no extension of the row is needed.
Under the growth conditions below, \(m_N\le2N=L_N\) for all sufficiently
large \(N\), and finitely many initial rows are irrelevant.

\begin{lemma}
\label{lem:moment-bound-dependent-condition}
Let \(p>2\). For a stationary row of length \(L_N\), condition
\eqref{eq:dependent-lindeberg-condition} follows from
\begin{equation}
\label{eq:p-moment-dependent-condition}
 \frac{L_Nm_N^{p-1}\E|X_{N,1}|^p}{s_N^p}\longrightarrow0.
\end{equation}
\end{lemma}

\begin{proof}
For \(x,a\ge0\),
\(x^2\1_{\{x>a\}}\le x^pa^{2-p}\). With
\(a=\varepsilon s_N/m_N\), stationarity gives
\[
\begin{split}
 &\frac{m_N}{s_N^2}\sum_{j=1}^{L_N}
 \E\!\left[X_{N,j}^2
 \1_{\{|X_{N,j}|>\varepsilon s_N/m_N\}}\right]\\
 &\qquad\le
 \varepsilon^{2-p}
 \frac{L_Nm_N^{p-1}\E|X_{N,1}|^p}{s_N^p}.
\end{split}
\]
The right-hand side tends to zero by
\eqref{eq:p-moment-dependent-condition}.
\end{proof}

\medskip
\noindent\textbf{Variance lower bound and moment verification.}

For fixed \(n\) and \(N\), write
\begin{equation}
\label{eq:finite-window-variance-definition}
 s_{n,N}^2(f):=\Var\bigl(Q_{n,N}(f)\bigr),
 \qquad s_{n,N}(f)\ge0.
\end{equation}

\begin{lemma}[Finite-window spectral variance bound]
\label{lem:finite-window-spectral-variance-lower}
For bounded periodic \(f\), \(n\ge0\), and \(N>nR\),
\begin{equation}
\label{eq:finite-window-spectral-variance-lower}
 s_{n,N}^2(f)\ge
 \lambda(2N-2nR-1)_+\rho_n(f)^2.
\end{equation}
If \(N\ge2Rn+1\), then
\(s_{n,N}^2(f)\ge\lambda N\rho_n(f)^2\).
\end{lemma}

\begin{proof}
Equations \eqref{eq:first-chaos-lower-bound-finiteN} and
\eqref{eq:Delta-equals-Pn-on-interior} give
\[
 s_{n,N}^2(f)\ge
 \lambda\int_{-N+nR}^{N-nR}|(\mathcal P_nf)(z)|^2\,\dd z.
\]
If \(g\ge0\) is \(1\)-periodic, its integral over an interval of length
\(L\) is at least
\((L-1)_+\int_0^1g\). Indeed, the interval contains a subinterval of
integer length \(\lfloor L\rfloor\), and the integral of a periodic
function over any interval of integer length equals that length times its
period average. Taking
\[
 g=|\mathcal P_nf|^2,\qquad L=2N-2nR,
\]
proves \eqref{eq:finite-window-spectral-variance-lower}. If
\(N\ge2Rn+1\), then
\[
 2N-2nR-1\ge N,
\]
which gives the last assertion.
\end{proof}

\begin{corollary}
\label{cor:spectral-growth-estimate-dependent}
Fix \(p>2\). If \(N\ge2Rn_N+1\) and \(\rho_{n_N}(f)>0\), then
\begin{equation}
\label{eq:spectral-growth-estimate-dependent}
 \frac{2Nm_N^{p-1}\E|X_{N,0}(f)|^p}{s_N(f)^p}
 \le C_{p,f}N^{1-p/2}
 \exp\{C_{p,f}D_N\log(D_N+2)\}\rho_{n_N}(f)^{-p}.
\end{equation}
In particular, for \(p=4\), the left-hand side tends to zero if
\begin{equation}
\label{eq:spectral-growth-condition}
 D_N\log(D_N+2)+
 \log_+\frac1{\rho_{n_N}(f)}=o(\log N).
\end{equation}
\end{corollary}

\begin{proof}
Proposition~\ref{prop:block-moment-estimate} gives
\[
 \E|X_{N,0}(f)|^p\le
 C_{p,f}^{p(D_N+1)}(1+D_N)^{p(D_N+1)}.
\]
Since \(n_N\le D_N\), \(m_N\le C_R(1+D_N)\). Lemma
\ref{lem:finite-window-spectral-variance-lower} gives
\(s_N(f)^p\ge\lambda^{p/2}N^{p/2}\rho_{n_N}(f)^p\).
Consequently,
\[
\begin{split}
 \frac{2Nm_N^{p-1}\E|X_{N,0}(f)|^p}{s_N(f)^p}
 &\le C_{p,f}N^{1-p/2}
 (1+D_N)^{p-1}\\
 &\quad\times
 C_{p,f}^{p(D_N+1)}
 (1+D_N)^{p(D_N+1)}
 \rho_{n_N}(f)^{-p}.
\end{split}
\]
All powers of \(1+D_N\), including the factor of degree \(p-1\), are
absorbed into the exponential factor
\[
 \exp\{C_{p,f}D_N\log(D_N+2)\}.
\]
This proves
\eqref{eq:spectral-growth-estimate-dependent}.

For \(p=4\), the logarithm of the right-hand side is bounded above by
\[
 -\log N+C_f+
 C_fD_N\log(D_N+2)+4\log_+\frac1{\rho_{n_N}(f)}.
\]
Under \eqref{eq:spectral-growth-condition}, it equals
\(-\log N+o(\log N)\) and hence tends to \(-\infty\).
\end{proof}

\medskip
% \noindent\textbf{A spectral increasing-depth central limit theorem.}

\begin{theorem}
\label{thm:growing-depth-clt-spectral}
Fix either recursion and a bounded real \(1\)-periodic \(f\). Suppose
\(\rho_{n_N}(f)>0\) eventually and
\begin{equation}
\label{eq:main-spectral-growth-condition}
 D_N\log(D_N+2)+
 \log_+\frac1{\rho_{n_N}(f)}=o(\log N).
\end{equation}
Then
\[
 \frac{Q_{n_N,N}(f)-\E Q_{n_N,N}(f)}
 {\sqrt{\Var(Q_{n_N,N}(f))}}
 \Longrightarrow\mathcal N(0,1).
\]
\end{theorem}

\begin{proof}
Condition \eqref{eq:main-spectral-growth-condition} implies
\[
 D_N\log(D_N+2)=o(\log N),
\]
so \(D_N=N^{o(1)}\), and in particular
\(n_N\le D_N=o(N)\). Hence, eventually,
\(N\ge2Rn_N+1\), \(m_N\le2N\), and
\(s_N(f)>0\) by Lemma
\ref{lem:finite-window-spectral-variance-lower}. By
\eqref{eq:growing-depth-block-sum},
\[
 Q_{n_N,N}(f)-\E Q_{n_N,N}(f)
 =\sum_{k=-N}^{N-1}X_{N,k}(f).
\]
Reindex with \(j=k+N+1\), so the row has length \(L_N=2N\).
Corollary~\ref{cor:spectral-growth-estimate-dependent} with \(p=4\),
followed by Lemma~\ref{lem:moment-bound-dependent-condition}, verifies
\eqref{eq:dependent-lindeberg-condition}. All hypotheses of
Theorem~\ref{thm:dependent-type-criterion} are now satisfied, and its
normalization \(s_N\) is exactly \(s_N(f)\).
\end{proof}

\medskip
\noindent\textbf{Proof of Theorem 1.4 and consequences.}

\begin{proof}[Proof of Theorem~\ref{thm:intro-growing-depth}]
Corollary~\ref{cor:zero-mode-first-chaos-nondegeneracy} gives
\(\rho_{n_N}(f)\ge\alpha_0^{D_N}|\widehat f(0)|\), and therefore
\[
 \log_+\frac1{\rho_{n_N}(f)}\le C_f(1+D_N).
\]
Condition \eqref{eq:intro-growing-depth-condition} implies
\eqref{eq:main-spectral-growth-condition}, so
Theorem~\ref{thm:growing-depth-clt-spectral} applies.
\end{proof}

In particular, the same conclusion holds for every bounded nonnegative
periodic \(f\) that is not zero almost everywhere.

\begin{corollary}
\label{cor:clt-one-nonvanishing-fourier-mode}
Suppose that \(\widehat f(k_0)\ne0\) for some \(k_0\in\Z\) and
\begin{equation}
\label{eq:one-mode-growth-condition}
 D_N\log(D_N+2)+
 \log_+\frac1{|p_{n_N,k_0}|}=o(\log N).
\end{equation}
Then the conclusion of Theorem
\ref{thm:growing-depth-clt-spectral} holds.
\end{corollary}

\begin{proof}
Equation \eqref{eq:rho-n-definition} gives
\(\rho_{n_N}(f)\ge|p_{n_N,k_0}||\widehat f(k_0)|\).
The fixed nonzero factor \(|\widehat f(k_0)|\) contributes only
\(O(1)\) after applying \(\log_+(1/\cdot)\). Thus
\eqref{eq:one-mode-growth-condition} implies
\eqref{eq:main-spectral-growth-condition}, and Theorem
\ref{thm:growing-depth-clt-spectral} applies.
\end{proof}

\begin{corollary}[Linear recursion]
\label{cor:linear-growing-depth-clt}
If
\[
 n_N\log(n_N+2)=o(\log N)
\]
and \(f\) has nonzero mean, then
\[
 \frac{Q_{n_N,N}^{\mathrm{lin}}(f)-\E Q_{n_N,N}^{\mathrm{lin}}(f)}
 {\sqrt{\Var(Q_{n_N,N}^{\mathrm{lin}}(f))}}
 \Longrightarrow\mathcal N(0,1).
\]
For a mean-zero \(f\), the same conclusion holds under
\eqref{eq:main-spectral-growth-condition} with \(D_N=n_N\).
\end{corollary}

\begin{proof}
For the linear recursion, \(D_N=d_{n_N}=n_N\). The nonzero-mean assertion follows from
Theorem~\ref{thm:intro-growing-depth}; the mean-zero assertion follows from
Theorem~\ref{thm:growing-depth-clt-spectral} with this value of \(D_N\).
\end{proof}

\begin{corollary}[Nonlinear recursion]
\label{cor:nonlinear-growing-depth-clt}
If
\[
 n_N2^{n_N}=o(\log N)
\]
and \(f\) has nonzero mean, then
\[
 \frac{Q_{n_N,N}^{\mathrm{nl}}(f)-\E Q_{n_N,N}^{\mathrm{nl}}(f)}
 {\sqrt{\Var(Q_{n_N,N}^{\mathrm{nl}}(f))}}
 \Longrightarrow\mathcal N(0,1).
\]
For a mean-zero \(f\), the same conclusion holds under
\eqref{eq:main-spectral-growth-condition} with \(D_N=2^{n_N}-1\).
\end{corollary}

\begin{proof}
Here
\[
 D_N\log(D_N+2)
 =(2^{n_N}-1)\log(2^{n_N}+1)\le Cn_N2^{n_N},
\]
so Theorem~\ref{thm:intro-growing-depth} applies.
\end{proof}

The explicit conditions can be read without the degree notation. In the
linear model, every sequence
\[
 n_N=o\!\left(\frac{\log N}{\log\log N}\right)
\]
satisfies the stated hypothesis. In the nonlinear model, if
\(\delta\in(0,1)\) is fixed and
\[
 n_N\le(1-\delta)\log_2\log N
\]
eventually, then
\[
 n_N2^{n_N}
 \le(1-\delta)(\log_2\log N)(\log N)^{1-\delta}
 =o(\log N).
\]
These examples illustrate only the sufficient scales delivered by the
proof.

\begin{remark}[Scope of the depth conditions]
The conditions in this subsection are sufficient. The local-count bound
\eqref{eq:block-moment-estimate} replaces the detailed recursive structure
by a power of one Poisson count, while the variance estimate retains only
the first chaos. The resulting gap between the moment upper bound and the
actual variance may be large. We establish neither necessity,
non-Gaussian behavior beyond these regimes, nor a rate of convergence, and
the conditions are not claimed to be optimal. In particular, the conclusions
do not identify a boundary between Gaussian and non-Gaussian behavior.
\end{remark}

\label{sec:examples}

\begin{example}[Depth zero]
At \(n=0\), \(W_0\equiv1\), and therefore
\[
 a_0=\lambda,\qquad b_0=0,\qquad
 \mathcal K_0=\lambda I.
\]
Thus \(\xi_0\) is white noise on \(\T\) with intensity \(\lambda\):
\[
 \E[\xi_0(f)\xi_0(g)]=\lambda\int_0^1f(x)g(x)\,\dd x.
\]
The insertion response is \(\Gamma_{0,z}=\delta_z\), whence
\(\gamma_0=\delta_0\), \(p_{0,k}=1\) for every \(k\), and
\[
 \theta_{0,k}=\lambda|p_{0,k}|^2=\lambda.
\]
Accordingly, the first-chaos bound is an equality at depth zero. This also
shows directly why the Gaussian limit is distribution-valued rather than
\(L^2\)-valued.
\end{example}

\begin{example}[Response at the first recursive depth]
The two recursions agree at depth one:
\[
 W_1(v)=\int_\R\varphi(u-v)\,\mu(\dd u).
\]
After inserting a point at \(z\),
\[
 W_1^{(z)}(v)-W_1(v)=\varphi(z-v).
\]
Consequently, the response measure in
\eqref{eq:Gamma-nz-definition} is
\[
\begin{split}
 \Gamma_{1,z}(\dd v)
 &=\varphi(z-v)\mu(\dd v)\\
 &\quad+
 \left\{\int_\R\varphi(u-z)\,\mu(\dd u)+\alpha_0\right\}
 \delta_z(\dd v).
\end{split}
\]
At \(z=0\), taking expectations by the Campbell formula gives
\begin{equation}
 \gamma_1(\dd s)
 =\lambda\varphi(-s)\,\dd s+
 \left(\lambda\int_\R\varphi(u)\,\dd u+\alpha_0\right)\delta_0(\dd s).
\end{equation}
If
\[
 \widehat\varphi(k):=\int_\R
 \varphi(u)e^{-2\pi iku}\,\dd u,
\]
then the response multiplier is
\begin{equation}
 p_{1,k}
 =\alpha_0+\lambda\int_\R\varphi(u)\,\dd u
  +\lambda\widehat\varphi(k),\qquad k\in\Z.
\end{equation}
This calculation also shows explicitly that a positive response measure may
have complex nonzero-frequency multipliers when the kernel is not symmetric;
the variance lower bound depends on \(|p_{1,k}|\).
\end{example}

\begin{example}[Constant mode]
For \(f=1\),
\[
 \Var(\xi_n(1))=\theta_{n,0}=\Sigma_n(1,1).
\]
Under Assumption~\ref{ass:local-positivity-kernel}, its logarithm is
\((2+o(1))n\log n\) in the linear model and
\((2\log2+o(1))n2^n\) in the nonlinear model.
\end{example}

\begin{example}[One nonzero frequency]
For \(k\ge1\), the Gaussian variables \(\xi_n(c_k)\) and
\(\xi_n(s_k)\) are independent and have common variance
\(\theta_{n,k}\). Their first-chaos contribution satisfies
\(\theta_{n,k}\ge\lambda|p_{n,k}|^2\). Thus a nonzero response
multiplier makes the corresponding frequency nondegenerate.
\end{example}

\begin{example}[Positive and mean-zero profiles]
If \(0<m_f\le f\le M_f\), then \(\Sigma_n(f,f)\) has the same leading
logarithmic exponent as \(\Sigma_n(1,1)\). If instead
\(\widehat f(0)=0\), the zero-mode lower bound is absent, and
\[
 \Sigma_n(f,f)\ge\lambda
 \sum_{k\ne0}|p_{n,k}|^2|\widehat f(k)|^2.
\]
The growing-depth problem is then controlled by the decay of the nonzero
multipliers visible to \(f\).
\end{example}

\appendix

\section{Proof of the uniform block moment estimate}
\label{sec:appendix-block-moments}

We first record an elementary moment estimate for a Poisson random variable.

\begin{lemma}[Poisson moment estimate]
\label{lem:poisson-moment-estimate}
If \(Z_a\sim\operatorname{Poisson}(a)\), where \(a\ge0\), then, for every
\(q\ge1\),
\begin{equation}
\label{eq:poisson-moment-estimate}
    \|Z_a\|_q\le q+a.
\end{equation}
\end{lemma}

\begin{proof}
The assertion is immediate when \(a=0\). Suppose that \(a>0\). For
\(t>0\) and \(x\ge0\),
\[
    x^q\le
    \left(\frac{q}{et}\right)^q e^{tx}.
\]
Since
\[
    \E[e^{tZ_a}]
    =\exp\!\bigl(a(e^t-1)\bigr),
\]
we obtain
\[
    \|Z_a\|_q
    \le
    \frac{q}{et}
    \exp\!\left(\frac{a(e^t-1)}q\right).
\]
Choose
\[
    t=\log\left(1+\frac qa\right).
\]
Then \(a(e^t-1)=q\), and hence
\[
    \|Z_a\|_q
    \le
    \frac{q}{\log(1+q/a)}.
\]
Finally,
\[
    \log(1+x)\ge\frac{x}{1+x},
    \qquad x\ge0,
\]
gives
\[
    \frac{q}{\log(1+q/a)}
    \le q+a.
\]
\end{proof}

\begin{proof}[Proof of Proposition~\ref{prop:block-moment-estimate}]
Fix \(p\ge1\). For \(v\in I_k\), Lemma
\ref{lem:pathwise-local-count} gives
\[
    W_n(v)
    \le
    M^{d_n}Z_{n,k}^{d_n},
    \qquad
    Z_{n,k}:=\mu([k-nR,k+1+nR]).
\]
Since \(\mu(I_k)\le Z_{n,k}\),
\[
\begin{aligned}
    |B_{n,k}(f)|
    &=
    \left|
      \int_{I_k}f(v)W_n(v)\,\mu(\dd v)
    \right| \\
    &\le
    \|f\|_\infty M^{d_n}Z_{n,k}^{d_n+1}.
\end{aligned}
\]
The random variable \(Z_{n,k}\) is Poisson distributed with mean
\[
    a_n^{\mathrm{loc}}
    :=\lambda(1+2nR).
\]
Consequently,
\[
\begin{aligned}
    \|B_{n,k}(f)\|_p
    &\le
    \|f\|_\infty M^{d_n}
    \|Z_{n,k}\|_{p(d_n+1)}^{d_n+1}.
\end{aligned}
\]
By Lemma~\ref{lem:poisson-moment-estimate},
\[
    \|Z_{n,k}\|_{p(d_n+1)}
    \le
    p(d_n+1)+\lambda(1+2nR).
\]
For both recursions, \(n\le d_n\). Therefore,
\[
    p(d_n+1)+\lambda(1+2nR)
    \le
    C_{p,\lambda,R}(1+d_n),
\]
where \(C_{p,\lambda,R}\) is independent of \(n\) and \(k\). It follows
that
\[
    \|B_{n,k}(f)\|_p
    \le
    C_{p,f}^{\,d_n+1}(1+d_n)^{d_n+1},
\]
where \(C_{p,f}\) depends only on
\[
    p,\lambda,R,M,\|f\|_\infty.
\]

By stationarity,
\[
\begin{aligned}
    \|X_{n,k}(f)\|_p
    &=
    \|B_{n,k}(f)-\E B_{n,k}(f)\|_p \\
    &\le
    \|B_{n,k}(f)\|_p
    +|\E B_{n,k}(f)| \\
    &\le
    2\|B_{n,k}(f)\|_p.
\end{aligned}
\]
Enlarging \(C_{p,f}\) absorbs the factor \(2\) and proves
\eqref{eq:block-moment-estimate}.
\end{proof}

% Optional acknowledgment (commented out in the original source):
% \section*{Acknowledgments}
% K. V. Hlyniana is the corresponding author and is also affiliated with the Institute of Mathematics, National Academy of Sciences of Ukraine.


\begin{thebibliography}{99}

\bibitem{BlaszczyszynYogeshwaranYukich2019}
Bart{\l}omiej B{\l}aszczyszyn, D.~Yogeshwaran, and Joseph~E. Yukich.
\newblock Limit theory for geometric statistics of point processes having
fast decay of correlations.
\newblock \emph{The Annals of Probability}, 47(2):835--895, 2019.
\newblock \href{https://doi.org/10.1214/18-AOP1273}
{doi:10.1214/18-AOP1273}.

\bibitem{Bogachev1998}
Vladimir~I. Bogachev.
\newblock \emph{Gaussian Measures}.
\newblock Mathematical Surveys and Monographs, vol.~62. American
Mathematical Society, Providence, 1998.
\newblock \href{https://doi.org/10.1090/surv/062}
{doi:10.1090/surv/062}.

\bibitem{DorogovtsevStochasticAnalysisHilbert}
A.~A. Dorogovtsev.
\newblock \emph{Stochastic Analysis and Random Maps in Hilbert Space}.
\newblock VSP, Utrecht, 1994.
\newblock ISBN 978-90-6764-163-0.

\bibitem{JansonPratelliRigo2024}
Svante Janson, Luca Pratelli, and Pietro Rigo.
\newblock Quantitative bounds in the central limit theorem for
\(m\)-dependent random variables.
\newblock \emph{ALEA, Latin American Journal of Probability and
Mathematical Statistics}, 21:245--265, 2024.
\newblock \href{https://doi.org/10.30757/ALEA.v21-10}
{doi:10.30757/ALEA.v21-10}.

\bibitem{LachiezeRey2019}
Rapha{\"e}l Lachi{\`e}ze-Rey.
\newblock Normal convergence of nonlocalised geometric functionals and
shot-noise excursions.
\newblock \emph{The Annals of Applied Probability}, 29(5):2613--2653,
2019.
\newblock \href{https://doi.org/10.1214/18-AAP1445}
{doi:10.1214/18-AAP1445}.

\bibitem{LachiezeReyPeccatiYang2022}
Rapha{\"e}l Lachi{\`e}ze-Rey, Giovanni Peccati, and Xiaochuan Yang.
\newblock Quantitative two-scale stabilization on the Poisson space.
\newblock \emph{The Annals of Applied Probability}, 32(4):3085--3145,
2022.
\newblock \href{https://doi.org/10.1214/21-AAP1768}
{doi:10.1214/21-AAP1768}.

\bibitem{LastPeccatiSchulte2016}
G{\"u}nter Last, Giovanni Peccati, and Matthias Schulte.
\newblock Normal approximation on Poisson spaces: Mehler's formula,
second order Poincar{\'e} inequalities and stabilization.
\newblock \emph{Probability Theory and Related Fields},
165(3--4):667--723, 2016.
\newblock \href{https://doi.org/10.1007/s00440-015-0643-7}
{doi:10.1007/s00440-015-0643-7}.

\bibitem{LastPenrose2011}
G{\"u}nter Last and Mathew~D. Penrose.
\newblock Poisson process Fock space representation, chaos expansion and
covariance inequalities.
\newblock \emph{Probability Theory and Related Fields},
150(3--4):663--690, 2011.
\newblock \href{https://doi.org/10.1007/s00440-010-0288-5}
{doi:10.1007/s00440-010-0288-5}.

\bibitem{LastPenrose2017}
G{\"u}nter Last and Mathew Penrose.
\newblock \emph{Lectures on the Poisson Process}.
\newblock Institute of Mathematical Statistics Textbooks, vol.~7.
Cambridge University Press, Cambridge, 2017.
\newblock \href{https://doi.org/10.1017/9781316104477}
{doi:10.1017/9781316104477}.

\bibitem{LastPenroseSchulteThaele2014}
G{\"u}nter Last, Mathew~D. Penrose, Matthias Schulte, and Christoph
Th{\"a}le.
\newblock Moments and central limit theorems for some multivariate Poisson
functionals.
\newblock \emph{Advances in Applied Probability}, 46(2):348--364, 2014.
\newblock \href{https://doi.org/10.1239/aap/1401369698}
{doi:10.1239/aap/1401369698}.

\bibitem{McLean2000StronglyElliptic}
William McLean.
\newblock \emph{Strongly Elliptic Systems and Boundary Integral Equations}.
\newblock Cambridge University Press, Cambridge, 2000.
\newblock ISBN 978-0-521-66375-5.

\bibitem{PeccatiSoleTaqquUtzet2010}
Giovanni Peccati, Josep Llu{\'i}s Sol{\'e}, Murad~S. Taqqu, and
Frederic Utzet.
\newblock Stein's method and normal approximation of Poisson functionals.
\newblock \emph{The Annals of Probability}, 38(2):443--478, 2010.
\newblock \href{https://doi.org/10.1214/09-AOP477}
{doi:10.1214/09-AOP477}.

\bibitem{PeccatiZheng2010}
Giovanni Peccati and Cengbo Zheng.
\newblock Multi-dimensional Gaussian fluctuations on the Poisson space.
\newblock \emph{Electronic Journal of Probability}, 15:1487--1527, 2010.
\newblock \href{https://doi.org/10.1214/EJP.v15-813}
{doi:10.1214/EJP.v15-813}.

\bibitem{Penrose2003}
Mathew~D. Penrose.
\newblock \emph{Random Geometric Graphs}.
\newblock Oxford Studies in Probability, vol.~5. Oxford University Press,
Oxford, 2003.
\newblock \href{https://doi.org/10.1093/acprof:oso/9780198506263.001.0001}
{doi:10.1093/acprof:oso/9780198506263.001.0001}.

\bibitem{Penrose2005}
Mathew~D. Penrose.
\newblock Multivariate spatial central limit theorems with applications to
percolation and spatial graphs.
\newblock \emph{The Annals of Probability}, 33(5):1945--1991, 2005.
\newblock \href{https://doi.org/10.1214/009117905000000206}
{doi:10.1214/009117905000000206}.

\bibitem{ReitznerSchulte2013}
Matthias Reitzner and Matthias Schulte.
\newblock Central limit theorems for \(U\)-statistics of Poisson point
processes.
\newblock \emph{The Annals of Probability}, 41(6):3879--3909, 2013.
\newblock \href{https://doi.org/10.1214/12-AOP817}
{doi:10.1214/12-AOP817}.

\bibitem{Shubin2001PDO}
Mikhail~A. Shubin.
\newblock \emph{Pseudodifferential Operators and Spectral Theory}.
\newblock 2nd ed. Springer, Berlin, 2001.
\newblock \href{https://doi.org/10.1007/978-3-642-56579-3}
{doi:10.1007/978-3-642-56579-3}.

\bibitem{Taylor2011PDEI}
Michael~E. Taylor.
\newblock \emph{Partial Differential Equations I: Basic Theory}.
\newblock 2nd ed., Applied Mathematical Sciences, vol.~115. Springer,
New York, 2011.
\newblock \href{https://doi.org/10.1007/978-1-4419-7055-8}
{doi:10.1007/978-1-4419-7055-8}.

\end{thebibliography}
\end{document}